\pdfoutput=1
\documentclass[11pt]{article}

\usepackage[T1]{fontenc}
\usepackage{lmodern}
\usepackage{microtype}
\usepackage{amsmath,amssymb,amsthm,mathtools}
\usepackage{xcolor}
\usepackage[margin=1.08in]{geometry}
\usepackage{enumitem}
\usepackage{graphicx}
\usepackage{float,needspace} 
\usepackage{booktabs,tabularx,array,longtable}
\usepackage{fancyvrb}
\usepackage[most]{tcolorbox}
\usepackage{tikz}
\usetikzlibrary{arrows.meta,positioning,calc,fit,shapes.geometric,decorations.pathreplacing}
\usepackage[
  backend=biber,
  style=alphabetic,
  sorting=nyt,
  giveninits=true,
  maxbibnames=99,
  maxalphanames=5,
  doi=true,
  isbn=false,
  url=false,
  eprint=true
]{biblatex}
\usepackage[
  colorlinks=true,
  linkcolor=magenta,
  citecolor=cyan,
  urlcolor=magenta,
  pdftitle={Notes on Generative AI for Mathematical Research},
  pdfauthor={Shurui Liu}
]{hyperref}
\usepackage[capitalise,noabbrev]{cleveref}

\setlist[enumerate]{leftmargin=2.2em,itemsep=2pt,topsep=4pt}
\setlist[itemize]{leftmargin=2em,itemsep=2pt,topsep=4pt}
\allowdisplaybreaks
\numberwithin{equation}{section}
\makeatletter
\renewcommand\@seccntformat[1]{%
  \ifcsname #1@cntformat\endcsname \csname #1@cntformat\endcsname
  \else \csname the#1\endcsname\quad \fi}
\newcommand\subsubsection@cntformat{\thesubsubsection.}
\newcommand\subsubsection@title[2]{
  \if\relax\detokenize{#2}\relax\else #1{\itshape #2\@addpunct{.}}\fi}
\newcommand\subsubsection@start{\@startsection{subsubsection}{3}{\z@}%
  {3.25ex \@plus1ex \@minus.2ex}{-.5em}{\normalfont\normalsize}}
\RenewDocumentCommand\subsubsection{s O{#3} >{\TrimSpaces}m}{%
  \IfBooleanTF{#1}
    {\subsubsection@start*{\subsubsection@title{}{#3}}}
    {\subsubsection@start[#2]{\subsubsection@title\enspace{#3}}}}
\makeatother

\theoremstyle{plain}
\newtheorem{theorem}{Theorem}[section]

\newtheorem{question}[theorem]{Question}
\newtheorem{conjecture}[theorem]{Conjecture}

\theoremstyle{definition}
\newtheorem{definition}[theorem]{Definition}
\newtheorem{principle}{Principle}

\theoremstyle{remark}
\newtheorem{remark}[theorem]{Remark}

\definecolor{primary}{HTML}{1F4E79}
\definecolor{accent}{HTML}{B5651D}
\definecolor{tint}{HTML}{EAF1F8}
\definecolor{warmtint}{HTML}{FBF1E6}
\definecolor{good}{HTML}{2E7D32}
\definecolor{bad}{HTML}{C62828}
\definecolor{muted}{HTML}{666666}

\newtcolorbox{takeaway}[1][Takeaway]{enhanced,breakable,colback=tint,
  colframe=primary,boxrule=0.6pt,arc=1.5mm,left=2mm,right=2mm,top=1mm,
  bottom=1mm,fonttitle=\bfseries\small,title={#1}}
\newtcolorbox{casebox}[1]{enhanced,breakable,colback=warmtint,
  colframe=accent,boxrule=0.6pt,arc=1.5mm,left=2mm,right=2mm,top=1mm,
  bottom=1mm,fonttitle=\bfseries\small,title={#1}}
\newtcolorbox{warnbox}[1][Caution]{enhanced,breakable,colback=bad!4,
  colframe=bad!70!black,boxrule=0.6pt,arc=1.5mm,left=2mm,right=2mm,top=1mm,
  bottom=1mm,fonttitle=\bfseries\small,title={#1}}
\newtcolorbox{promptbox}{enhanced,breakable,colback=white,colframe=muted,
  boxrule=0.4pt,arc=0pt,left=2mm,right=2mm,top=1mm,bottom=1mm,
  fontupper=\small\ttfamily\raggedright}
\makeatletter
\AddToHook{env/takeaway/before}{\if@noskipsec\leavevmode\par\fi}
\AddToHook{env/casebox/before}{\if@noskipsec\leavevmode\par\fi}
\AddToHook{env/warnbox/before}{\if@noskipsec\leavevmode\par\fi}
\AddToHook{env/promptbox/before}{\if@noskipsec\leavevmode\par\fi}
\makeatother

\newcommand{\XX}{\mathbb{X}}
\newcommand{\sys}[1]{\textsc{#1}}
\newcommand{\file}[1]{\texttt{#1}}

\newcommand{\yes}{{\color{good}\ensuremath{\checkmark}}}
\newcommand{\no}{{\color{bad}\ensuremath{\times}}}
\newcommand{\TNR}{\mathrm{TNR}}
\newcommand{\FNR}{\mathrm{FNR}}

\DefineVerbatimEnvironment{shell}{Verbatim}{fontsize=\small,frame=single,
  framesep=2mm,rulecolor=\color{muted}}

\tikzset{
  box/.style={draw=black!35,fill=white,rounded corners=1.2mm,line width=0.5pt,
    inner sep=2mm,align=center,font=\small},
  pbox/.style={box,draw=primary,line width=0.7pt},
  abox/.style={box,draw=accent,line width=0.7pt},
  flow/.style={-{Stealth[length=2mm]},line width=0.7pt,color=black!55},
  pflow/.style={flow,color=primary},
  aflow/.style={flow,color=accent},
  lab/.style={font=\footnotesize,color=muted},
}

\title{Notes on Generative AI for Mathematical Research}
\author{%
  Shurui Liu\\
  \href{mailto:srliu@stanford.edu}{\texttt{srliu@stanford.edu}}\\[1em]
  {\normalsize\textit{Department of Mathematics, Stanford University}}\\[0.4em]
  {\normalsize November 2026}
}
\date{}

\begin{document}
\maketitle

\begin{abstract}
These notes describe how generative AI is now used to solve and to check
research-level mathematics, with a focus on agents and their harnesses. For the proving agents \sys{Rethlas}, \sys{Danus}, and
\sys{Pharos} and the verification agent \sys{OSVerify}, we explain the design
principles, compare each with its base model, present case studies, and show
how to use and contribute to them. We also briefly review the history of AI
for mathematics and the community's reactions to it.

\medskip\noindent\textit{2020 Mathematics Subject Classification.}
Primary 68T42; Secondary 68T20, 68T50, 68V15, 00A35, 01A61.
\end{abstract}

\tableofcontents
\clearpage 

\section{Introduction}\label{sec:introduction}
\subsection{Overview}
\subsubsection{Motivation}
These notes summarize my exploration of AI for mathematical research in 2026.
They are expanded from talks I gave at Stanford University, Duke Kunshan
University, and Cursor, to audiences of both mathematicians and computer
scientists, and they also draw on my experience teaching Math~289, a course on
AI and mathematics at Stanford. Their purpose is to record my understanding
and to help bridge the information gap between mathematicians and AI
researchers.

\subsubsection{Natural-language agents for mathematical research}
We focus on the \textit{informal} side of AI for mathematics rather than the
formal side. That is, problem solving, proof generation, and verification
here mean reasoning in \textit{natural language} (such as English), not in a
formal language (such as Lean~4).

\subsubsection{Bias and limitations}
Most of the material and the views in these notes come from my own research
in mathematics and in AI for mathematics up to November 2026. They are
therefore \textit{biased} and represent only my own perspective. They are
\emph{not} meant to represent the mathematical community, nor to cover all
related work or approaches.
\subsubsection{Agents discussed in this paper}
\Cref{tab:systems} gives an overview of the systems discussed in
\cref{sec:systems,sec:verification}. All of them are open source.
\Cref{sec:systems} treats the systems that solve problems and prove theorems:
Rethlas, Danus, and Pharos. \Cref{sec:verification} treats OSVerify, a
verification agent that checks proofs. For each system we first compare its
capabilities with those of the base model, then describe the design, then
discuss case studies, and finally its limitations.

\begin{table}[ht]
\centering\small
\renewcommand{\arraystretch}{1.2}
\begin{tabularx}{\textwidth}{@{}l>{\raggedright\arraybackslash}X>{\raggedright\arraybackslash}X l@{}}
\toprule
\textbf{System} & \textbf{Idea} & \textbf{What it mimics} & \textbf{Reference} \\
\midrule
Rethlas & a generator and a verifier in a loop & one mathematician and a referee & \cite{ju2026rethlas} \\
Danus & a research group of agents sharing a DAG of verified facts & a research group of junior students with a administrator and an external expert consultant & \cite{liu2026danus} \\
Pharos & re-designed Danus for stronger models: route board, numerical tools, deeper thinking, etc. & a research group of top PhDs with a PI & \cite{dong2026jets} \\
OSVerify & one agent proposes errors, another validates them & a careful and strict referee & \cite{nie2026osverify} \\
\bottomrule
\end{tabularx}
\caption{The systems discussed in \cref{sec:proving,sec:verification}.}
\label{tab:systems}
\end{table}
\Needspace{0.5\textheight}
\subsubsection{Selected comparison}
\Cref{tab:results-selected} collects representative case studies comparing base models with our systems. More details are discussed in \S\ref{sec:systems} and \S\ref{sec:verification}.
\begin{table}[H]
\centering\small
\renewcommand{\arraystretch}{1.25}
\begin{tabularx}{\textwidth}{@{}>{\raggedright\arraybackslash}X>{\raggedright\arraybackslash}p{22mm}>{\raggedright\arraybackslash}p{36mm}>{\raggedright\arraybackslash}p{42mm}@{}}
\toprule
\textbf{Problem} & \textbf{Field} & \textbf{Model alone} & \textbf{With our system} \\
\midrule
Cocharacters \cite{ju2026rethlas} & algebraic groups & GPT-5.5 Pro: {\color{bad}wrong proof} & Rethlas (GPT-5.5): {\color{good}correct proof, 44\,min} \\
Matroid KL polynomials \cite{cheng2026kl} & combinatorics & GPT-5.6 Sol Ultra: {\color{bad}incomplete} & Rethlas (GPT-5.6 Sol): {\color{good}counterexample, 5\,h} \\
Matroid tangent classes \cite{cheng2026tangent} & combinatorics & GPT-5.5 Pro: {\color{bad}no solution} & Danus (GPT-5.5): {\color{good}solved, $\sim$5 days} \\
Hartshorne conjecture, reproved \cite{dong2026jets} & algebraic geometry & GPT-6-Astra with expert guidance: {\color{bad}no solution} & Pharos (GPT-6 Astra): {\color{good}new proof via jets, 20\,h, autonomous} \\
Manin constant, review \cite{chi2026manin} & number theory & Claude Opus~5.5, Claude Fable~5.1, Codex (GPT-6 Astra): {\color{bad}12 rounds, no error found} & OSVerify (GPT-5.6 Sol): {\color{good}one error found, fixes suggested} \\
11 OpenAI preprints \cite{openai2026math} & various & OpenAI's internal model: {\color{bad}mistakes overlooked} & OSVerify (GPT-6-Astra): {\color{good}found mistakes in each} \\
\bottomrule
\end{tabularx}
\caption{Base models alone versus our systems: selected results.}
\label{tab:results-selected}
\end{table}
\subsection{Organization of the paper}
\Cref{sec:history} briefly reviews the history of AI for mathematical
research, reactions of the math community, and some philosophical
perspectives. \Cref{sec:research} discusses what
mathematical research is and why we believe it differs from competition
mathematics and from chess. \Cref{sec:agents} reviews the basic notions of
LLMs, harnesses, and agents, and discusses the motivation and design
principles of agents for mathematics. \Cref{sec:systems,sec:verification}
describe our agents for problem solving and for verification: their designs,
comparisons with the base models, and case studies. \Cref{sec:use} explains
how to use and contribute to these systems, and \cref{sec:future} points out
some future directions.
\subsection{Declaration of human usage and AI usage}
This paper was written primarily by the author. The tables and figures were generated with AI tools under the author's direction and then revised by the author. Generative AI tools were also used to improve the exposition (correcting typos and grammar, rephrasing sentences), to detect imprecisions and to help find relevant references. All philosophical views and conceptual understanding presented here were developed by the author through independent thought and discussions with other people; some of the technical computer science background was learned with the help of AI tools. The author takes full responsibility for the correctness of the content.

\subsection{Acknowledgements}
These notes are expanded from talks at the AI for Mathematics seminar and in
Math~289 at the Stanford Department of Mathematics, at the Graphics Cafe
seminar in the Stanford Department of Computer Science, at the AI for
Mathematics seminar at Duke Kunshan University, and at Cursor. S.\,Liu thanks
the organizers and hosts, and the audiences for inspiring discussions,
especially Otis Chodosh, Andrea Montanari, Ravi Vakil, Lexing Ying, George
Nakayama, Zander Majercik, Jiajun Ma, Xiaoqian Xu, and Shengtong Zhang.

S.\,Liu thanks his collaborators for teaching him through practice, especially Bin Dong, Guoxiong Gao, Haocheng Ju, James Zou, Haotian Ye, Fan Nie, Bin Wu, Zeming Sun, and Jiedong Jiang. He also thanks Brian Conrad, Fred Rajasekaran, Henry Bosch, Ronnie Cheng, Haosen Wu, Rosie Shen, Vaughan McDonald, Richard Taylor, Yanan Jiang, Travor Liu, Qiyue Tang, Zhuo Zhang, Carl Schildkraut, Michael Ren, Zeyu Wang, Romain Speciel, Judson Kuhrman, Selim Amar, Yizhen Chen, Hongrui Chen,  and Alec Shelley, for meaningful discussions, and his advisor, Xinwen Zhu, for constant support.

\section{A brief history of AI for mathematics}\label{sec:history}
This section is not meant to be a comprehensive survey of AI-assisted results. We select some examples, and discuss what they indicate and our personal sentiment about the rapid progress.

\subsection{The old era: computer-assisted proofs}
Long before large language models and agents, there were many
computer-assisted proofs.
\subsubsection{Examples}
\paragraph{Exhaustive symbolic search in equational logic}
In 1996, McCune's prover EQP showed that every Robbins algebra is Boolean, a
question open since the 1930s \cite{mccune1997robbins}.
\paragraph{Exhausting a huge finite number of cases}
Appel and Haken's proof of the Four Color Theorem (1976) \cite{appel1977four}
started a long debate on whether a proof that no human can survey is a proof.
The computer, however, relied on a human proof strategy that reduced the
original problem to finitely many cases.

\paragraph{Formal languages and interactive theorem provers}
Lean~\cite{demoura2015lean} is an interactive theorem prover: definitions, statements, and proofs are written in a formal language and mechanically checked by a small trusted kernel. Its community-maintained library Mathlib~\cite{mathlib2020} provides a large, unified body of formalized mathematics, covering much of the undergraduate curriculum and substantial graduate-level material. The Liquid Tensor Experiment, posed as a challenge by
Scholze~\cite{scholze2022liquid} and completed by a team led by Commelin, was to formalize in Lean
a difficult theorem of Clausen and Scholze on liquid vector spaces.
\subsubsection{Brute force and deterministic programs}
All the examples above are, more or less, brute-force searches by
deterministic programs. The computer did not propose new ideas; it only
computed faster than humans. Given enough time, a human could in principle
carry out the same task.
\subsection{The not-so-old era: machine learning methods}

From about 2016, machine learning entered mathematics in several roles.

\subsubsection{Guiding Lean provers} Neural premise selection \cite{irving2016deepmath}
and then language models trained on formal proofs \cite{polu2020gptf} learned
to suggest the next tactic or lemma. Tools such as LeanDojo
\cite{yang2023leandojo} and semantic search engines for Mathlib
\cite{gao2024semantic} were built with machine learning. A sequence
of increasingly strong formal provers followed (DeepSeek-Prover, Kimina, Seed-Prover \cite{chen2025seed}, Aristotle \cite{achim2025aristotle}), as did
statement autoformalization \cite{gao2024herald,wang2025aria}.

\subsubsection{Numerical discovery of examples} Davies et al.\ \cite{davies2021advancing} used
attribution methods on trained networks to guide mathematicians to a new
theorem in knot theory and to progress on the combinatorial invariance
conjecture for Kazhdan--Lusztig polynomials. Wagner used reinforcement
learning to find counterexamples in graph theory \cite{wagner2021constructions}.
PatternBoost
\cite{charton2024patternboost} and AlphaEvolve \cite{novikov2025alphaevolve}
extended the idea to many constructions in combinatorics and analysis. In these systems the model proposes and a program \emph{scores}; they work best when the objective is a computable number.

\subsubsection{Symbolic deduction for IMO geometry} AlphaGeometry \cite{trinh2024alphageometry}
combined a language model with a symbolic deduction engine and solved 25 of 30
IMO geometry problems. In July 2024 AlphaProof and AlphaGeometry~2 reached
silver-medal level at the IMO, with proofs in Lean \cite{hubert2025olympiad}.
\subsubsection{The boundary between intelligence and brute force}
These methods go beyond the deterministic programs of the previous era: a learned model, rather than a fixed rule, decides which tactic to try, which example to examine, or which construction to add next. This is a form of intuition, and it makes searches feasible that brute force could not finish. Yet the examples above still share features with brute force. The search takes place in a narrow space fixed in advance; every candidate is checked automatically, by the Lean kernel, a scoring program, or a symbolic engine. The model is to prune a search tree. The examples of the next era differ in nature although where exactly the line between clever search and genuine intelligence lies is debatable.
\subsection{The language-model and agent era}\label{sec:llm-era}

\subsubsection{Math competitions and low-hanging Erdős problems}
\begin{itemize}
  \item \textbf{July 2025: IMO gold.} An experimental general-purpose
  reasoning model from OpenAI solved five of the six problems at the IMO
  (35/42, gold-medal level), writing natural-language proofs under contest
  time limits; the proofs were graded by former IMO medalists rather than by
  the IMO \cite{openai2025imo}. Google DeepMind's Gemini Deep Think obtained
  the same score with official certification.
  \item \textbf{January 2026: Erdős problems.} GPT-5.2 Pro produced proofs of
  several Erdős problems, formalized in Lean by Harmonic's Aristotle
  \cite{achim2025aristotle}; Tao's public wiki tracks AI contributions to
  Erdős problems \cite{tao2026erdoswiki}. These are generally regarded as
  low-hanging fruit. Some had in fact already been solved in the literature:
  Erdős Problem~\#397, for instance, turned out to be Problem~3 of Day~1 of
  the 2012 Chinese IMO Team Selection Test. Others were open but minor and
  yielded to known methods, such as Erdős Problem~\#728.
\end{itemize}

\subsubsection{Combinatorics and small research results}
\begin{itemize}
  \item \textbf{February 2026: Aletheia.} Google DeepMind's Aletheia, a
  research agent built on Gemini Deep Think that generates, verifies, and
  revises proofs in natural language, produced a paper on eigenweights in
  arithmetic geometry without human intervention \cite{feng2026eigenweightsarithmetichirzebruchproportionality}, contributed to several
  papers written with mathematicians, and autonomously solved four open
  Erdős problems in a study of 700 of them \cite{feng2026towards}.
  \item \textbf{February 2026: First Proof.} Eleven mathematicians released
  10 lemmas from their current research \cite{abouzaid2026firstproof}.
  OpenAI's internal model produced solutions judged correct for about
  5 of them after expert review \cite{openai_first_proof_2026}, and
  Aletheia had 6 accepted by a majority of assessors
  \cite{feng2026aletheia}.
  \item \textbf{August 2026: a Hadamard matrix of order 668.} On 12 August,
  L.~Alpöge, P.~Voinov and S.~Reynolds-Haertle, assisted by Claude, announced
  Hadamard matrices of all twelve orders below 2000 for which none was known
  \cite{alpoge2026hadamard}. These include order 668, the smallest open case
  of the Hadamard conjecture since the construction of order 428 in 2005.
  The structure of the matrices was later analyzed by Egan, Kotsireas,
  London and Ó~Catháin \cite{egan2026hadamard}.
\end{itemize}

\subsubsection{Publishable and important results}
Several online lists keep track of major AI-assisted results in mathematical
research, for example \cite{mathoverflow502120}. Here we list only some that
received public attention.
\begin{itemize}
  \item \textbf{May 2026: unit distances.} An internal OpenAI model disproved
  Erdős's conjecture that $n$ points in the plane determine at most
  $n^{1+o(1)}$ unit distances, by a construction from algebraic number theory
  (number fields and a Golod--Shafarevich argument)
  \cite{openai2026unitdistance}. A group of nine mathematicians verified and
  explained the argument \cite{alon2026remarks}. Because of the novelty of the
  method, many mathematicians consider the result worthy of a top journal
  such as the \emph{Annals}. Days later, Anthropic reported that Claude
  Mythos had independently found a (slightly weaker) disproof.
  \item \textbf{July 2026: the Jacobian conjecture.} L.~Alpöge (Anthropic)
  announced an explicit polynomial map $\mathbb C^3 \to \mathbb C^3$ of
  degree~7 with constant Jacobian determinant that is not injective, found
  with Claude Fable~5, disproving the Jacobian conjecture in every dimension
  $n \ge 3$. The counterexample was verified in Isabelle within a day
  \cite{ramos2026jacobian} and explained by Tao \cite{tao2026jacobian}. The
  conjecture is famous, but experts mostly care about the two-dimensional
  case, which remains wide open.
  \item \textbf{August 2026: ten advances.} On 1 August, OpenAI released
  results of its unreleased model Astra on ten open problems in mathematics
  and theoretical computer science, with Lean certificates for all ten
  \cite{openai2026ten}. They include the existence of a non-sofic group and a
  disproof of Connes's rigidity conjecture. However, researchers in the field raised the possibility that key ideas came from their unpublished work.
  \item \textbf{August 2026: a complex structure on $S^6$.} On 23 August,
  L.~Alpöge posted a 108-page manuscript, produced by an internal version of
  Claude under his direction, constructing a one-parameter family of compact
  complex threefolds diffeomorphic to $S^6$. This answers Hopf's 1947
  question of whether the six-sphere admits a complex structure. P.~Engel
  gave a self-contained exposition \cite{engel2026s6}, and B.~Alexeev
  formalized the proof in Lean. As of this writing, experts consider the
  construction plausible, and scrutiny continues.
  \item \textbf{September 2026: Navier--Stokes with a force term.} OpenAI
  released a proof of finite-time blowup for the \emph{forced}
  three-dimensional Navier--Stokes equations, with a Lean certificate of about
  230{,}000 lines, produced by a large swarm of agents
  \cite{openai2026navierstokes}. It is not the unforced Millennium Problem: the
  forced version is expected to be much easier, and it is not as interesting
  as the unforced one. It was also disputed on credit and on whether the
  formal statement faithfully captures the informal one
  \cite{bastounis2026lost}.
  \item \textbf{October 2026: 722 manuscripts.} OpenAI released 722
  manuscripts in 372 families of results, produced by an internal model
  from about 4{,}000 attempted problems, credited to ``OpenAI'' with no human
  authors; about half were unverified at release \cite{openai2026math}. Within
  days, three were withdrawn for a sign error and fourteen were revised
  \cite{retractionwatch2026openai}. The claimed results include full
  dimension for Kakeya sets in $\mathbb R^4$ and a quasi-Riemann hypothesis.
  If correct, many of them would be publishable in top journals and would be
  milestones in their areas.
  \item \textbf{November 2026: the Riemann hypothesis.} On 20 November,
  L.~Alpöge announced that, working with Claude, he had found a proof of the
  Riemann hypothesis, and that the method extends to the generalized Riemann
  hypothesis \cite{alpoge2026remark,sciam2026riemann}. The proof was released
  not as a conventional paper but as a Lean formalization, reportedly of
  about 17 million lines, together with a short explanatory note.
  \item \textbf{November 2026: Navier--Stokes without a force term.} The same
  day, T.~Buckmaster and L.~Alpöge posted a heavily AI-assisted proof of
  finite-time blowup for the \emph{unforced} three-dimensional Navier--Stokes
  equations, with a Lean library of about 600{,}000 lines
  \cite{buckmaster2026navierstokes}. As of this writing it has not been
  independently checked, and the numerical part of the proof has not been
  released.
\end{itemize}

\subsection{Sentiment about AI for mathematics}\label{sec:sentiment}

\subsubsection{A divided community}
Opinions in the mathematical community are sharply divided. For example, see the communal blog \emph{Proofs and
Prompts} \cite{proofsandprompts} and also the AMS pages on mathematics and artificial intelligence \cite{ams2026ai}. For the views of graduate students and postdocs, see the survey \cite{brandenberger2026earlycareer}. Reactions span a full spectrum. At one extreme is total opposition. For example, in August 2026, M.~Weinreich made the case for total opposition to the use of AI in mathematics \cite{weinreich2026crisis}. The Association for Human Mathematics, founded in September 2026, aims to protect mathematics as a human endeavor \cite{ahmath}. At the other extreme is wholesale promotion of AI-generated mathematics without human involvement. In between, many people call for using AI with transparency and responsibility; see, for example, the recommendations in \cite{commelin2026shaping}. As base models improve and more important problems are solved, these opinions keep evolving.

\subsubsection{An analogy from science fiction}
Cixin Liu's \emph{Remembrance of Earth's Past}
trilogy~\cite{liu2014threebody,liu2015darkforest} begins when humanity learns
that the Trisolarans, a far more advanced alien civilization, have launched a
fleet that will reach Earth in roughly four centuries. Ahead of the fleet they
send ``sophons'', particle-sized computers that observe everything humans say
or write, though not what they think. Most of the story concerns the long
wait: how a society reorganizes itself around a superior intelligence that has
not yet arrived. People split into factions that agree on the facts and
disagree on what follows. Adventists want the aliens to replace a humanity
they regard as beyond repair, Redemptionists revere the aliens and hope to
save their world, triumphalists and defeatists argue over whether resistance
can succeed, and a few ``Wallfacers'' are asked to plan in the one place that
cannot be observed, their own minds. The second volume adds the image of the
universe as a ``dark forest'', in which any civilization that reveals itself
risks destruction.

We find the parallel with mathematics useful for three reasons. First, the
situation has the same shape: a non-human intelligence whose abilities are
growing on an uncertain timetable, and a community that must decide how to
respond before it knows how capable that intelligence will become. Second,
the interesting variation is among the human responses, which is what a
taxonomy of factions records. Third, the sophons have a direct counterpart:
anything a mathematician posts can be read, and possibly used, by AI systems
and used in their training. However, the analogy is to the pattern of reactions only. The labels below are caricatures that describe attitudes without endorsing or condemning them. We borrow the names loosely, and take a few from elsewhere where the novels run out.

\subsubsection{A taxonomy of reactions}
\begin{description}
  \item[Adventists] hold that AI has already superseded human mathematics
    entirely, so that doing mathematics by hand is a waste of time and
    resources. A milder version holds that human mathematics is itself flawed,
    or at least too slow~\cite{darpa2025expmath}, and that AI working alongside
    humans will accelerate it enormously and so rescue us.

  \item[Redemptionists] believe we should spend our effort improving AI, and
    help it escape the ``convex hull'' of its training data and the
    ``discrimination'' against AI-generated work in publications.

  \item[Triumphalists] expect a quick victory~\cite{liu2015darkforest}. They
    believe that AI's current abilities are a passing phase: it will hit a
    ceiling sooner or later, or it can never leave the convex hull of its
    training data, or industry will lose interest and stop investing. In each
    case humans regain the lead.

  \item[Capitulationists] concede that human mathematics still has value and a
    future, but see none in themselves, and conclude that all they can do is
    embrace AI to make a living.

  \item[The Protracted-War school] accepts that AI keeps improving, but
    observes that it is still weak at exposition, theory-building and real
    innovation. Its strategy is to keep finding what AI does badly and to add
    human value there.

  \item[Wallfacers]~\cite{liu2015darkforest} keep their thoughts to
    themselves, like their namesakes in the novel, who plan only inside their
    own heads because everything spoken or written is overheard. For them
    mathematics has become a dark forest: any idea that is disclosed may be
    taken and scooped, by other people or by AI, so ideas must stay secret
    until the work is finished.

  \item[Wanderers,] named after the \emph{xiaoyao} (``free and easy
    wandering'') of the \emph{Zhuangzi}~\cite{zhuangzi2013watson}, simply enjoy
    life and the pleasure of understanding mathematics
    (cf.~\cite{thurston2006proof}). Their model is literature, which has
    professional writers and amateurs alike, and whose golden age comes only
    when everyone can read and write.

  \item[The Steam faction] sees an industrial revolution in mathematics, in
    which the subject enters an era of mass production, close to what Tao has
    described as mathematics on an industrial scale~\cite{wong2024uncharted}.
\end{description}

\subsubsection{Evaluation and hiring}\label{sec:hiring}
Whatever one's position, the community's systems of evaluation and hiring
will have to change. Some hiring committees have begun to publish advice and describe changes in how they assess applications \cite{proofsandprompts2026hiring}. For doctoral training, a group of symplectic geometers has published recommendations on the use of AI for PhD students and supervisors \cite{proofsandprompts2026phd}.

\subsection{Philosophy perspectives}
\subsubsection{Some general opinions}
The mathematical community holds a range of opinions on using AI in research.
I suspect that most people would agree with the following.
\begin{enumerate}
    \item A diversity of ideas and perspectives is good for mathematics and
    for humanity.
    \item AI is already very good at solving problems and generating proofs,
    although its ability varies across problems and areas. But mathematical
    research is not only about problem solving.
    \item AI is not yet a magic box that solves every problem and completes
    every task, although people predict different futures.
    \item Mathematicians value intellectually interesting and challenging
    work; work that machines can easily do will eventually lose its value
    (say, at least when hiring committees evaluate candidates; see
    \cref{sec:hiring}).
    \item One should always at least intend to serve the advancement of
    mathematics, science, and humanity, even though people see this
    differently and act differently.
    \item Honesty and reputation always matter in the community.
\end{enumerate}

\subsubsection{Some personal opinions}
My own position consists of the list above together with the following.
\begin{enumerate}
    \item The advancement of mathematics and the survival of mathematicians
    may diverge. We should certainly aim at the former. But mathematicians,
    especially those who prefer not to use AI and choose to rely on their own
    creativity, should be respected and protected, since no one can be sure
    that they will not be needed in the long run.
    \item Humans, or at least some humans, can still add value on top of AI.
    For this to remain true, we must act responsibly toward the younger
    generation, who are shaped by education, by the community's system of
    evaluation, by social media, and by leading mathematicians.
    \item In an AI-assisted proof, three groups of people contribute to the
    human understanding of the mathematics: (i) those who build the AI tools,
    (ii) those who solve the problem with the tools, and (iii) those who
    understand the solution and help others understand it. I believe their
    contributions compare as (i) $\ll$ (ii) $<$ (iii).
\end{enumerate}

\subsubsection{Suggestions}
These notes focus on technical practice rather than philosophical debate. Still,
I would like to make a few gentle suggestions. I believe these mistakes are
usually made without bad intent, since AI tools make them easy to make.
\begin{itemize}
    \item Please take care not to scoop other people's work. An AI system
    will not ask whether others are working on the same problem, so we have
    to ask ourselves. In particular, if someone has announced a result and
    its main strategy in a seminar or at a conference, or has shared their
    ideas with you, please do not use AI to finish the work before they do.
    \item Please avoid producing piles of AI slop in areas you do not
    understand, and please do not advertise or hype a result before it has
    been carefully verified. An AI system writes just as fluently when it is
    wrong, and in an unfamiliar area it is hard to tell the difference.
    \item If an AI system solves a weak or uninteresting variant of a famous
    conjecture, please do not advertise or hype it as major progress on the
    conjecture itself. An AI system will not tell you whether a variant is
    interesting; that is for experts in the area to judge.
\end{itemize}

\section{What is mathematical research}\label{sec:research}

This section is non-technical. In \cref{subsec:research-loop} we decompose mathematical research into six stages, and in \cref{subsec:competitions} we discuss why research-level mathematics differs from competitions.
\subsection{The research loop}\label{subsec:research-loop}
\subsubsection{The six stages}
\begin{figure}[ht]
\centering
\begin{tikzpicture}[
    stage/.style={box,draw=black,minimum width=19mm,minimum height=12mm,inner xsep=1mm,font=\footnotesize\bfseries,text=black},
    hstage/.style={stage,draw=primary,fill=tint,line width=1.1pt,text=black,font=\footnotesize\bfseries},
    note/.style={font=\footnotesize,color=black,align=center,text width=24mm}]
  \node[stage] (disc) at (0,0) {Discovery};
  \node[hstage,right=3.8mm of disc] (solve) {Problem\\solving};
  \node[hstage,right=3.8mm of solve] (ver) {Verification};
  \node[stage,right=3.8mm of ver] (und) {Understanding};
  \node[stage,right=3.8mm of und] (comm) {Communication};
  \node[stage,right=3.8mm of comm] (app) {Application};
  \foreach \name/\desc in {
      disc/{propose\\conjectures},
      solve/{(dis)prove\\or compute},
      ver/check the solution,
      und/{conceptually\\internalize},
      comm/write and teach,
      app/{joy, job, \ldots}}
    \node[note,below=1.5mm of \name] {\desc};
  \foreach \a/\b in {disc/solve,solve/ver,ver/und,und/comm,comm/app}
    \draw[flow,color=black] (\a) -- (\b);
  \draw[flow,color=black,rounded corners] (app.north) -- ++(0,6mm) -| (disc.north);
\end{tikzpicture}
\caption{The research loop. The systems in this paper target the two
highlighted stages, problem solving and verification, which have
traditionally been the most visible part of research.}
\label{fig:loop}
\end{figure}
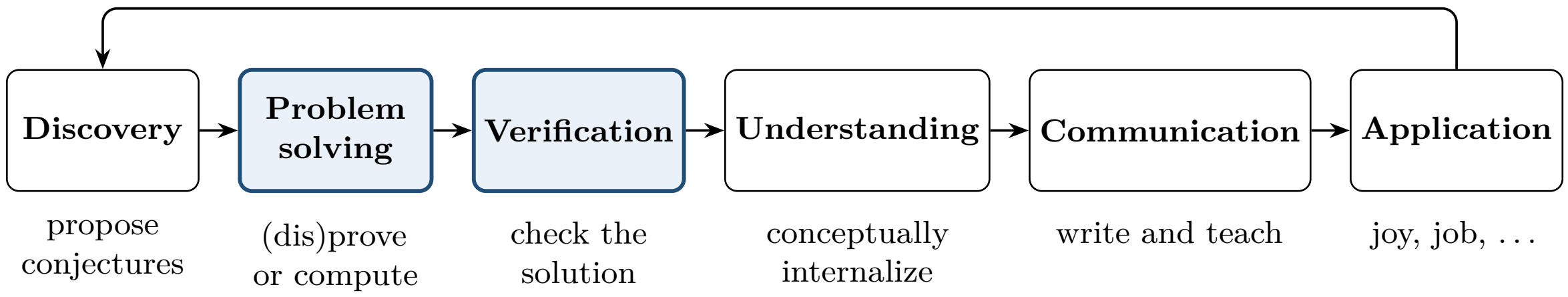

Mathematical research can be roughly decomposed into a loop of six stages
(\cref{fig:loop}).

\begin{enumerate}
  \item \textbf{Discovery:} Observing phenomena and proposing conjectures to explain them in a unified or conceptual way.
  \item \textbf{Problem solving:} Proving or disproving a conjecture, or computing an answer.
  \item \textbf{Verification:} Checking whether a solution is correct.
  \item \textbf{Understanding:} Internalizing why a proof works, which idea is
  essential and novel, and how it connects to other results. 
  \item \textbf{Communication:} Writing papers, giving talks, and teaching.
  \item \textbf{Application:} Applications to other problems within mathematics and to other subjects. This stage also includes the joy that mathematical exploration brings.
\end{enumerate}
\subsubsection{Other expository writings}
Thurston's essay \emph{On proof and progress in mathematics} (1994) \cite{thurston2006proof} is a beautiful account of what understanding means for mathematicians. I also personally like Larry Guth's essay \emph{What math means to me} \cite{guth2026mathmeans}. 
\subsubsection{Beyond problem solving}
The systems in these notes target problem solving and verification. Research, however, is much more than problem solving: a good mathematician is not the same as a good problem solver.

\subsection{Why research is different from competitions}\label{subsec:competitions}
\subsubsection{The chess analogy and competition math}
Mathematics is sometimes compared with chess. The comparison fits competitions (IMO, Putnam, \ldots) reasonably well: the syllabus is fixed; every problem is guaranteed to have a solution; and there is a large archive of past problems to train on. 
\subsubsection{How research differs}
However, research problems differ from competition problems in at least five ways.
\begin{enumerate}
\item \textbf{Open scope.} There is no fixed set of moves, and no guarantee that a
  solution exists, or that it can be reached with methods from the existing literature. A statement may even be undecidable from the standard axioms, as the continuum hypothesis is in ZFC.
\item \textbf{Generalization.} In competitions, the techniques behind past problems are reused on new ones. In research, a technique that solved one problem often does not transfer, and a new problem may need a brand new idea.
\item \textbf{Long horizon.} Proofs run from a few pages to hundreds, and
  explorations take days for an agent and months or years for a human. A wrong lemma or a bad route wastes
  everything built on it.
\item \textbf{Sparse reward.} It is hard to tell whether a partial result or step
  is meaningful, or to confirm that a long proof is fully correct.
\item \textbf{Evolving and noisy data.} Papers on arXiv and in journals contain errors, and they record finished proofs, not the reasoning or intuition that led to them or the attempts that failed. Moreover, the
  frontier is by definition not in the training data, so a system must learn at test time.
\end{enumerate}
\subsubsection{Math research varies across areas}
Research also varies greatly from one area of mathematics to another. For example:
\begin{itemize}
	\item Combinatorics tends to be technique-heavy, without a single central theme.
	\item The geometric and categorical Langlands programs tend to be theory-heavy, organized around a central theme and a guiding dream.
	\item Low-dimensional topology and symplectic geometry tend to be intuition-heavy; some geometric arguments are hard to write down rigorously, and visualization is essential.
\end{itemize}

It is unclear whether ability in one area, for humans or for models, transfers to another. Yet the areas interact: an insight or result from one field may solve a problem in another.
\subsubsection{Transfer to other domains}
For computer scientists, these properties make research mathematics a
valuable testbed for long-horizon agents and reasoning: the tasks are hard, yet correctness is objective in principle and the domain is free of physical-world noise. Techniques that work here are likely to transfer to other long-horizon tasks, such as automated research in software engineering, GPU-kernel optimization, and theoretical physics.

\section{What are harnesses and agents}\label{sec:agents}
Roughly speaking, harness engineering is the art of using large language models (LLMs) to solve specific tasks.
\subsection{What is an LLM}
In one sentence, an LLM is a machine that predicts the next token, where a token is roughly a word. We first recall what a token is.
\subsubsection{What is a token}
The word \textit{token} has different meanings in different contexts. In natural language processing (NLP) and large language models (LLMs), a token is a basic unit of text, typically a word or a piece of a word. For example:
\begin{itemize}
    \item \textit{apple}.
    \item \textit{un}, \textit{believ}, \textit{able}.
\end{itemize}
There are also special tokens, such as \texttt{<EOS>}, which marks the end of a sentence.\par 
\subsubsection{Tokens are determined before model training}
The vocabulary $\mathcal{V}$ is the set of all tokens, and it is learned from a text corpus before training and then fixed once and for all. Together with the vocabulary one obtains a deterministic program, the \emph{tokenizer}, that translates between text and tokens. One of the most popular choices is the byte-pair encoding (BPE) tokenizer; see \cite{raschka2025implementingabytepairencodingbpetoke} for details. 
\subsubsection{Three views of an LLM}
There are several ways to think about an LLM:
\begin{itemize}
    \item An LLM is a machine that completes a sentence by predicting the next token to say. 
    \item An LLM is a function
    \[
        p_{\theta}\colon \mathcal{V}^{\textrm{Fin}} \longrightarrow \Delta(\mathcal{V}),
    \]
    where $\mathcal{V}^{\textrm{Fin}}$ is the set of finite sequences of tokens, and $\Delta(\mathcal{V})$ is the set of probability distributions over $\mathcal{V}$. It completes a sentence \emph{autoregressively}: it repeatedly samples the next token given the input tokens and the tokens generated so far.
    \item An LLM $p_{\theta}$ is a map $\mathcal{V}^{\textrm{Fin}} \longrightarrow \Delta(\mathcal{V})$ built by composing fixed nonlinear functions with linear maps, where $\theta$ denotes the parameters of the linear maps. Training updates $\theta$ so that the model's predictions improve.
\end{itemize}
\begin{remark}
\textit{Large language model} is a fixed term: ``large'' does not necessarily mean that the model, that is, its parameter set $\theta$, is large.
\end{remark}
\subsection{A first definition: agent = model + harness}
There are also several ways to define agents and harnesses. The following is the most popular.
\subsubsection{The definition from computer science}
\begin{definition}
Agent = model + harness.
\end{definition}
In the words of V.~Trivedy (LangChain, 2026): ``If you're not the model,
you're the harness'' \cite{trivedy2026harness}.

\begin{remark}
    The defining equation involves two undefined terms. It reminds me of Newton's second law $F=ma$ (force and mass), and both are useful in practice. This one splits the work of improving an AI system into two parts: improving the model (pre-,
mid-, and post-training to update $\theta$), and improving the harness ($\theta$ unchanged).
\end{remark}

\subsubsection{The harness lets the model use tools}
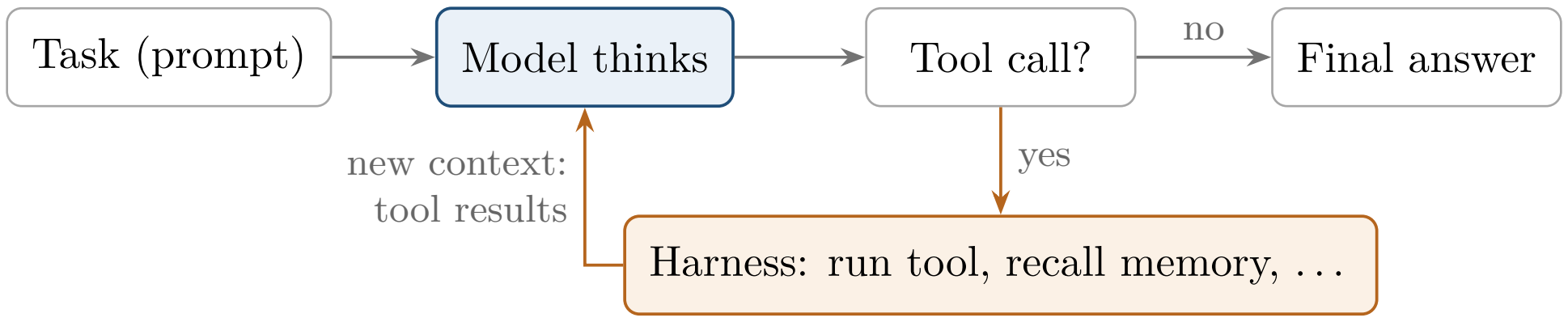
\begin{figure}[ht]
\centering
\begin{tikzpicture}[x=1cm,y=1cm]
  \node[box,minimum width=22mm,minimum height=8mm] (u) at (0,0) {Task (prompt)};
  \node[pbox,fill=tint,minimum width=24mm,minimum height=8mm] (m) at (3.4,0) {Model thinks};
  \node[box,minimum width=22mm,minimum height=8mm] (q) at (6.8,0) {Tool call?};
  \node[box,minimum width=22mm,minimum height=8mm] (a) at (10.2,0) {Final answer};
  \node[abox,fill=warmtint,minimum width=46mm,minimum height=8mm] (t) at (6.8,-1.7) {Harness: run tool, recall memory, \ldots};
  \draw[flow] (u) -- (m);
  \draw[flow] (m) -- (q);
  \draw[flow] (q) -- node[above,lab]{no} (a);
  \draw[aflow] (q) -- node[right,lab]{yes} (t);
  \draw[aflow] (t.west) -| node[pos=0.75,left,lab,align=right]{new context:\\tool results} (m.south);
\end{tikzpicture}
\caption{The tool-use loop. The model sees only its context, and the harness decides what enters the context next.}
\label{fig:toolloop}
\end{figure}

As in \cref{fig:toolloop}, the model reads its context and produces text. When the output contains a structured \emph{tool call}, say
\begin{Verbatim}[fontsize=\small,frame=single,framesep=2mm,rulecolor=\color{muted}]
<tool_call>{"name": "search", "query": "quasi-complete local ring"}</tool_call>
\end{Verbatim}
the harness executes the call (here a literature search), appends the result
to the context, and the model continues generating tokens.

\subsection{A second definition: harness = context management}
\subsubsection{Definition from first principles}
A language model is, at bottom, a machine for continuing text. Everything it
does depends on its context window. This suggests a working definition from first principles.

\begin{definition}
A harness is a \emph{context-management} system: a system external to the model whose aim is to give the model, at every step, exactly the context it needs for the task, namely everything useful and nothing irrelevant or misleading.
\end{definition}
\subsubsection{A ``set-theoretic'' definition of harness}
A harness may contain the following components.
\begin{description}
  \item[Infrastructure:] a sandbox and a workspace in which the agent can
  read and write files and run programs safely; process management (run, session and job IDs; CPU, memory, disk and wall-clock limits per process); browser and computer use.
  \item[Loop support:] keeping the agent going until a goal is met;
  compaction of the context.
  \item[Tools:] code execution, Lean, web search, and so on, usually exposed through
  the Model Context Protocol (MCP).
  \item[Memory:] plain text (Markdown files), append-only structured notes (searchable JSONL files), or SQLite databases, organized in short- and long-term tiers.
  \item[Prompts and skills:] role descriptions (\file{AGENTS.md},
  \file{CLAUDE.md}) and \emph{skills}, which are Markdown instruction files
  that give guidance or procedures for a task.
  \item[Orchestration:] how several agents collaborate, and how work is
  distributed among them.
\end{description}

\subsubsection{All components serve context management}
\begin{figure}[ht]
\centering
\begin{tikzpicture}[x=1cm,y=1cm,
    comp/.style={box,text width=44mm,minimum height=13mm,inner sep=1.2mm,font=\small}]
  \fill[tint!70,rounded corners=2mm] (-2.7,0.95) rectangle (13.1,-2.65);
  \foreach \name/\x/\y/\head/\sub in {
      in/0/0/Infrastructure/{sandbox, workspace},
      lp/5.2/0/Loop support/{compaction, keep going},
      tl/10.4/0/Tools/{search, code execution, MCP},
      me/0/-1.7/Memory/{md, JSONL, sqlite},
      sk/5.2/-1.7/Prompts and skills/{roles, \texttt{SKILL.md} files},
      or/10.4/-1.7/Orchestration/{subagents, handoff, routing}}
    \node[comp] (\name) at (\x,\y) {\textbf{\head}\\[-0.1ex]{\footnotesize\color{muted}\sub}};
  \draw[-{Stealth[length=2.5mm]},line width=1.3pt,color=primary] (5.2,-2.7) -- node[right,lab]{all serve} (5.2,-3.3);
  \node[pbox,fill=tint,minimum width=60mm,minimum height=8mm,font=\small\bfseries] at (5.2,-3.75) {Context management};
\end{tikzpicture}
\caption{Components of a harness.}
\label{fig:components}
\end{figure}
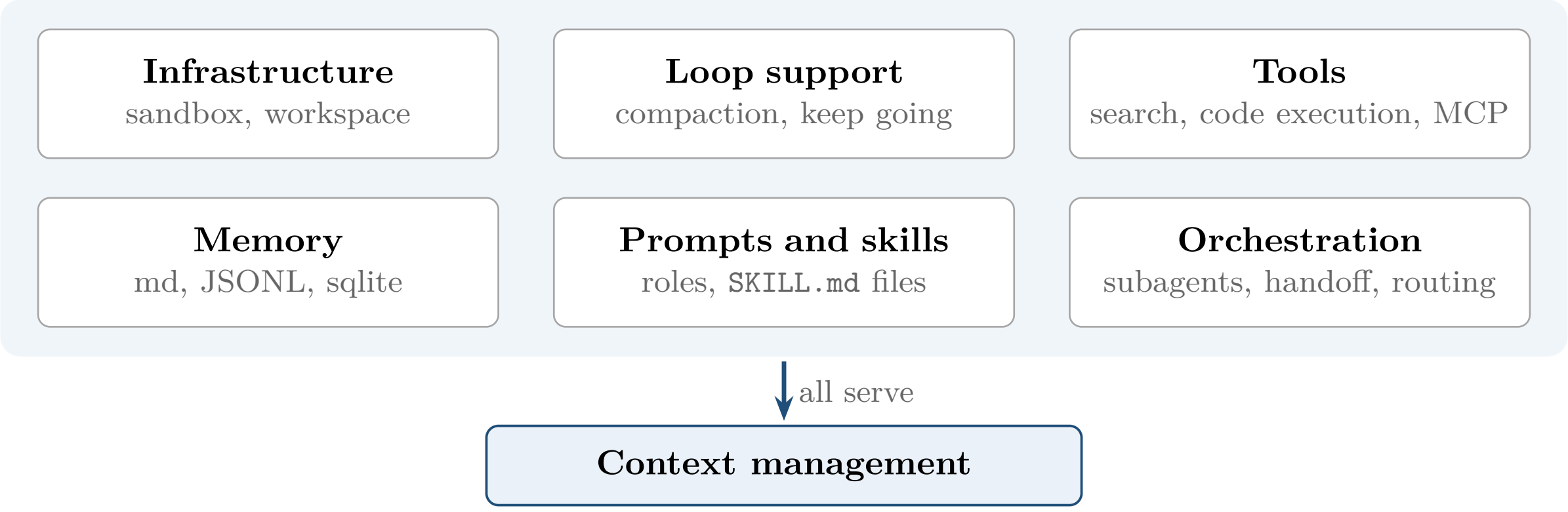

All the components of a harness (\cref{fig:components}) ultimately serve
context management.

\subsection{General-purpose harnesses (Codex, Claude Code) and their limitations}\label{sec:general-harness}
Most readers have already used a harness. 
\subsubsection{Frontier models ship inside general-purpose harnesses}
Claude Fable~5, GPT-6 Astra and other frontier models all ship inside
general-purpose harnesses:
\begin{itemize}
  \item \textbf{Chat interfaces} (ChatGPT Pro, Claude): extended thinking,
  web search, code execution, and memory managed in the cloud.
  \item \textbf{Coding agents} (Codex, Claude Code): they run shell commands and read and write files on your own computer, with memory files kept there too.
\end{itemize}

\subsubsection{Built-in controls are small harness features}
\begin{itemize}
    \item \texttt{/goal}: keeps the agent working until a goal is reached.
    \item Reasoning-effort levels up to \texttt{max}: deeper reasoning per step.
    \item \texttt{ultra} mode: runs more subagents in parallel.
\end{itemize}

\subsubsection{What still needs improvement}\label{sec:what-goes-wrong}

Even the strongest model in a strong general harness still needs improvement in at least four places.

\begin{enumerate}
  \item \textbf{Workflow for math research.} Their multi-agent orchestration
  and procedures are built for general purposes, not for research mathematics.
  \item \textbf{Context engineering.} Irrelevant or misleading information
  should not enter the context. 
  \item \textbf{Long horizon.} With several lines of work in one context
  (topic A, then topic B), the model loses track of which one it is
  pursuing.
  \item \textbf{Efficient use of compute.} A better-engineered harness can use compute more effectively than \texttt{/goal} or \texttt{ultra}, and can align better with human standards, saving both money and human effort.
\end{enumerate}

\subsubsection{Summary}
The harness matters at least whenever the task exceeds the context window,
or when the context confuses the model's attention.

\subsection{Why work on harnesses and agents}
\subsubsection{Why harnesses rather than model training}
\begin{enumerate}
  \item \textbf{The CS reason.} For most tasks,
  training-free methods should be tried before training. A harness is cheaper
  and faster to iterate on than a model. Moreover, an effective harness can help with training, by generating trajectories to use as synthetic training data, or through agentic reinforcement learning.
  \item \textbf{The financial reason.} Training a model with frontier mathematical ability is not a viable path for most research groups because of its huge cost. An effective harness scales compute more efficiently and gets more out of a model on the same budget. It may also make cheaper or open-source models capable of certain tasks.
  \item \textbf{Other reasons.} A good harness relies on an understanding of the base model's strengths and weaknesses, and designing one produces intuitions about how base models should improve. It may also improve our high-level understanding of the tasks themselves.
\end{enumerate}
\subsubsection{The harness evolves as foundation models improve}
As models improve, they absorb the simpler parts of the harness and make
more ambitious ones possible. A step-by-step skill instruction written for a weaker model may
be unnecessary for a stronger one, which can instead exploit a harness built
for harder tasks, such as recursive self-improvement.
\subsubsection{Which components survive stronger models?}
Components that \emph{compensate for weaknesses} of the model (such as detailed
skills or retrieval tricks) tend to become unnecessary as models improve; Pharos, for example, removed several of Danus's. Components that \emph{encode structure} or \emph{manage context} (such as bounded context per agent, or separation of roles) do not obviously depend on model strength. Whenever a task is longer than the context window, some harness will be needed.

\subsection{Design principles: how to design harnesses and agents}\label{sec:principles}

\subsubsection{Beyond trial and error}
Designing an effective harness is largely empirical and requires trial and
error. Still, a few principles can guide the design. The systems of
\cref{sec:proving,sec:verification} were designed with the three below.
\subsubsection{First principles: context engineering}
\begin{principle}\label{pr:context}
Provide the correct context: everything the model needs to solve the task, and nothing irrelevant.
\end{principle}
For example, a validation agent should start from a fresh context that does not contain the reasoning trace of the solution; and low-quality information should be compacted away.
\subsubsection{Mathematical intuition: mimic the human workflow}
\begin{principle}\label{pr:mimic}
Mimic how human experts solve the task.
\end{principle}
For example, human habits are distilled into \emph{skills}, external resources become \emph{tools}, and the notebook becomes \emph{memory}. The multi-agent orchestration of Danus and Pharos mimics a research group.
\subsubsection{Computer-science intuition: scale the compute meaningfully}

\begin{principle}\label{pr:scale}
Scale compute without too much waste, and compute can buy intelligence.
\end{principle}

The hope is that 10 smart people with IQ 140, organized into a well-designed collaboration, can match one person with IQ 400.


\makeatletter
\let\origsubsubsection\subsubsection
\let\origparagraph\paragraph
\renewcommand\subsubsection@cntformat{\thesubsubsection\quad}
\renewcommand\subsubsection{\@startsection{subsubsection}{3}{\z@}%
  {-3.25ex\@plus -1ex \@minus -.2ex}{1.5ex \@plus .2ex}{\normalfont\normalsize\bfseries}}
\providecommand\paragraph@cntformat{}
\renewcommand\paragraph@cntformat{\theparagraph.}
\newcommand\paragraph@start{\@startsection{paragraph}{4}{\z@}%
  {3.25ex \@plus1ex \@minus.2ex}{-.5em}{\normalfont\normalsize}}
\RenewDocumentCommand\paragraph{s O{#3} >{\TrimSpaces}m}{%
  \IfBooleanTF{#1}
    {\paragraph@start*{\subsubsection@title{}{#3}}}
    {\paragraph@start[#2]{\subsubsection@title\enspace{#3}}}}
\makeatother
\setcounter{secnumdepth}{4}
\addtocontents{toc}{\protect\setcounter{tocdepth}{3}}
\section{Agents for proving and problem solving: Rethlas, Danus, Pharos}\label{sec:systems}\label{sec:proving}
This section discusses three agents for proving or disproving conjectures: Rethlas, Danus, and Pharos. For each agent we begin with results comparing the base model alone with the harness, then describe the design, examine case studies in detail, and conclude with the limitations that motivate the next system.

\subsection{Rethlas: a generate--verify loop}\label{sec:rethlas}
Rethlas is fully open source, with a technical report \cite{ju2026rethlas} and code at \url{https://github.com/frenzymath/Rethlas}.

\subsubsection{Comparison with base models}
\paragraph{Comparison}
\Cref{tab:rethlas-results} compares the same base model with and without the
Rethlas harness on three problems. In each case Rethlas was given only the statement.

\begin{table}[ht]
\centering\small
\renewcommand{\arraystretch}{1.25}
\begin{tabularx}{\textwidth}{@{}>{\raggedright\arraybackslash}X>{\raggedright\arraybackslash}p{24mm}>{\raggedright\arraybackslash}p{38mm}>{\raggedright\arraybackslash}p{38mm}@{}}
\toprule
\textbf{Problem} & \textbf{Field} & \textbf{Model alone} & \textbf{With Rethlas} \\
\midrule
Rational conjugacy of cocharacters \cite[\S5.1.1]{ju2026rethlas} & algebraic groups & GPT-5.5 Pro: {\color{bad}wrong proof} & GPT-5.5: {\color{good}correct proof, 44\,min} \\
Unimodality of matroid KL polynomials \cite{cheng2026kl} & combinatorics & GPT-5.6 Sol Ultra: {\color{bad}incomplete} & GPT-5.6 Sol: {\color{good}counterexample, 4\,h\,59\,min} \\
Manin constant conjecture \cite{chi2026manin} & algebraic number theory & GPT-6 Astra: {\color{bad}needed many rounds of human guidance} & GPT-6 Astra, with OSVerify as verifier: {\color{good}correct proof, fully automated} \\
\bottomrule
\end{tabularx}
\caption{The same model with and without the harness. ``Model alone'' means
the model's web interface at its strongest setting.}
\label{tab:rethlas-results}
\end{table}
\paragraph{Other results}
A partial, growing list of results obtained with Rethlas is maintained at
\url{https://github.com/frenzymath/Rethlas_results}. For example:
\begin{itemize}
  \item \textit{commutative algebra}: Anderson's question \cite[\S4]{ju2026rethlas}; eight
  problems from the Cahen--Fontana--Frisch--Glaz collection
  \cite{cahen2014open,jiang2026openproblems}; questions of Erman--Sam on
  Boij--Söderberg theory;
  \item \textit{$p$-adic Hodge theory}: a refined lift-independence theorem \cite{pan2026lift};
  \item \textit{analysis and probability}: two problems from Brezis's list
  \cite{dou2026brezis}; injectivity of phase retrieval;
  \item \textit{combinatorics}: inverse Kazhdan–Lusztig polynomials of matroids need not be log-concave; inverse Z-polynomials of loopless matroids have unimodal coefficients.
  \item \textit{relative geometric Langlands} (with OSVerify as verifier): the automorphic commutator relation
  \cite{liu2026automorphic}; properties of Plancherel algebras
  \cite{hexagon2026plancherel}.
\end{itemize}

\subsubsection{Agentic workflow design}

\begin{figure}[ht]
\centering
\includegraphics[width=\linewidth]{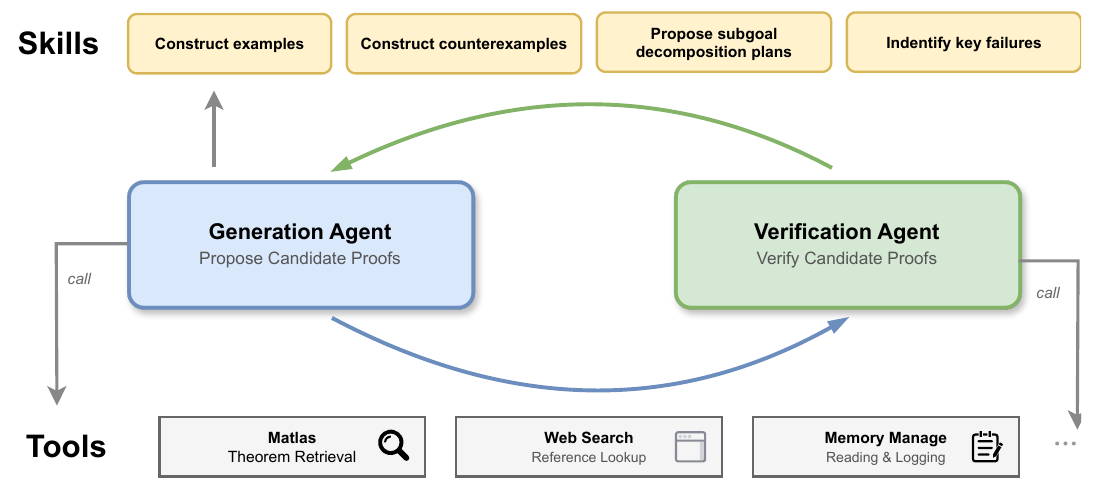}
\caption{The workflow design of Rethlas \cite{ju2026rethlas}.}
\label{fig:rethlas}
\end{figure}
Rethlas (\cref{fig:rethlas}) consists of two agents.
\paragraph{Generation} The \textbf{generation agent} writes a candidate proof, called a
  \emph{blueprint}: a paper-like Markdown file of lemmas, each with a statement and a proof, ending with the main theorem in its verbatim original form.
\paragraph{Verification} The \textbf{verification agent} referees the blueprint. Each submission is judged by an agent with a \emph{fresh} context, which returns a structured report: a summary, a list of critical errors, a list of gaps, and repair hints.
\paragraph{Iteration} The generate-verify loop repeats until acceptance, which requires zero critical
  errors and zero gaps.

\subsubsection{Design principles in Rethlas}
The design follows the principles of \cref{sec:principles} closely.

\paragraph{Mimic the human workflow: skills are habits} Rethlas encodes some habits of a mathematician attacking a hard problem as \emph{skills}: Markdown instruction files that the generation agent loads when it judges them relevant.
\begin{center}\small
\begin{tabular}{@{}ll@{}}
\toprule
\textbf{Exploration} &  construct toy examples; construct counterexamples; \\
& search relevant results; query memory \\
\textbf{Planning} & propose subgoal decomposition plans \\
\textbf{Proving} & obtain immediate conclusions; direct proving; recursive proving \\
\textbf{Iteration} & identify key failures; verify proof \\
\bottomrule
\end{tabular}
\end{center}
There is no fixed order: the agent loops over \emph{assess the state, choose a
skill, act}. In \emph{recursive proving}, competing decomposition plans run
in parallel as subagents, each told the difficulties found in all plans; if
all fail, the \emph{identify key failures} skill extracts their common
obstructions to seed the next round of plans. 

Here is an example of a skill:
\begin{Verbatim}[fontsize=\footnotesize,frame=single,framesep=2mm,rulecolor=\color{muted},label={construct-toy-examples}]
## Procedure
1. Construct simpler cases (low degree, small dimension,
   special forms, canonical objects).
2. Ensure the toy example satisfies all assumptions of
   the target statement or subgoal.
3. Check that the conclusion also holds in the toy example.
\end{Verbatim}
The verifier has its own skills. It checks \emph{the smallest deduction steps} for validity, applicability, and missing hypotheses, and it checks that \emph{references are correct}: external results must be cited with their full statement and source.
\paragraph{Mimic the human workflow: search the literature, and think independently} The central tool is \textbf{Matlas}
(\url{https://matlas.ai}), a semantic search engine over mathematical
statements. It is built from 8.51 million statements drawn from 435K published papers in 180 journals and over 1{,}900
textbooks, with dependency extraction so that retrieved statements are self-contained \cite{ju2026rethlas}. Rethlas blocks search in odd iterations and re-enables it in even ones, to balance \emph{thinking deeply} against \emph{learning from the literature}.

\paragraph{Mimic the human workflow: the notebook is memory} Raw logs are scratch paper. Structured
notes are kept in append-only JSONL channels per problem
(\file{toy\_examples}, \file{counterexamples}, \file{failed\_paths},
\file{verification\_reports}, \ldots) and retrieved by BM25 search. 

\paragraph{Mimic the human workflow: the author is not the referee} Authors believe their own proofs,
and so do models; self-checking suppresses errors only mildly. The
mathematical community's answer is the referee, someone with no stake in the proof being true.
\paragraph{Context engineering: continuous author, fresh referee}
\begin{itemize}
  \item the \emph{generator} has a \emph{continuous} context, so intuition, partial progress, and dead ends persist across iterations;
  \item the \emph{verifier} has a \emph{fresh} context per submission: it only sees the blueprint and the references, and nothing of how the blueprint was produced, so it is not distracted by the reasoning trace.
\end{itemize}
\paragraph{Scale compute effectively} The generation agent searches over candidate solutions. Because failed paths are recorded, compute is not wasted on repeating old mistakes.
\subsubsection{Case studies: copyright and understanding}
The technical report \cite{ju2026rethlas} contains three case studies: an open problem of D.~D.~Anderson in commutative algebra, a problem in $p$-adic Hodge theory \cite{pan2026lift}, and the algebraic-group problem discussed below. The report analyzes the reasoning trajectories of Rethlas; here we discuss instead two phenomena that we observed.
\paragraph{Case study 1: hallucination or galaxy-brain move?}
The prompt is the following algebraic group question.
\begin{question}[suggested by B.~Conrad]
Let $G$ be a connected smooth affine group over a field $k$, and
$\XX_*(G) = \operatorname{Hom}_{k\text{-gp}}(\mathrm{GL}_1, G)$. For a field
extension $K/k$, is the map $\XX_*(G)/G(k) \to \XX_*(G_K)/G(K)$ injective?
\end{question}

The answer was known for reductive $G$ and for
perfect $k$; the general case was not in the literature. Both Rethlas and GPT-5.5 Pro solved the easy cases. On the general case, GPT-5.5 Pro (web) gave a wrong proof with no meaningful progress, while Rethlas (GPT-5.5 \texttt{xhigh}) proved injectivity in 44 minutes.

The strategy was to reduce to the quasi-reductive case (the agent
coined the phrase ``wound unipotent case'') and then comparing relative and absolute roots and Weyl groups.

There was, however, one small mistake. The proof needs a key ingredient, \cite[Theorem C.2.15]{cgp2015}, from a copyrighted book. Rethlas requires the agent to quote the full statement of every external result it cites, which here would have meant reproducing copyrighted text. So the agent did not cite \cite{cgp2015}. Instead it searched the web for \texttt{CGP Thm C.2.15}, found the survey \cite{conradprasad2017} of \cite{cgp2015} on Conrad's homepage, and cited \cite[Theorem 5.3.2(iv)]{conradprasad2017}. The proof of that theorem refers back to \cite[Theorem C.2.15]{cgp2015}, but its statement is partial and does not contain the assertion the proof needed.

\paragraph{Case study 2: can the agent extract the method?}

This case study shows that understanding still needs human experts: the agent solved exactly the task it was given, but the general method behind its solution had to be extracted by humans.

There are two conjectures in matroid theory:
\begin{conjecture}[Gedeon--Proudfoot--Young 2017]\label{conjreal}
The Kazhdan--Lusztig polynomial $P_M(t)$ of a matroid $M$ is real-rooted.
\end{conjecture}
\begin{conjecture}[Elias--Proudfoot--Wakefield 2016]\label{conjlogcon}
The Kazhdan--Lusztig polynomial $P_M(t)$ of a matroid $M$ is log-concave.
\end{conjecture}

We tried only the first one with Rethlas. The prompt to Rethlas contained only the definition of $P_M$ and the instruction
``prove or disprove that all zeros of $P_M(t)$ lie on the negative real
axis.'' Rethlas (GPT-5.6 Sol, \texttt{max}) found an explicit counterexample to Conjecture~\ref{conjreal} in 4\,h\,59\,min. GPT-5.6 Sol Ultra in the web interface, given two prompts, did not finish.

However, the agent did only the minimum needed for the task it was given. The human authors of \cite{cheng2026kl} turned its construction into a method (\cite[Theorem 1.1]{cheng2026kl}), gave it a
geometric interpretation (\cite[\S3]{cheng2026kl}), and used it to prove that $P_M$ need not even be \emph{unimodal}, which refutes Conjecture~\ref{conjreal}, \ref{conjlogcon} at once. The agent did not
generalize its construction, explain it geometrically, or notice that it bore on the second conjecture.

\subsubsection{Why a generate--verify loop is not enough}

A single loop has several limitations, which motivate the next systems.
\paragraph{No decomposition of the problem} A mathematician
  accumulates many small facts and adjusts strategy as they come in. Rethlas pursues one hypothesis at a time, as a single depth-first search with no mixture of expertise.
\paragraph{Linear logical structure} A problem may need several people
  on different routes, and the routes interact, so the logical structure of the work need not be linear. A single Markdown blueprint cannot hold this, and two agents cannot safely edit it at once, so the loop cannot simply be parallelized.
\paragraph{Context engineering} Every round, the generator rereads the whole draft, with all verified partial results \emph{and their proofs}, even when only the statements are needed and the lemmas should be used as black boxes.
\paragraph{Scale compute to buy intelligence} Rethlas has two components, generation and verification, and each can be scaled. Danus (\cref{sec:danus}) and Pharos (\cref{sec:pharos}) scale compute on generation, from a single author--referee loop to a research group. OSVerify (\cref{sec:osverify}) scales compute on verification, to obtain a stricter verifier.


\subsection{Danus: a research group of agents}\label{sec:danus}

Danus is fully open source, with a paper \cite{liu2026danus} and code at \url{https://github.com/frenzymath/Danus}. Since GPT-6 Astra came out, Danus has been superseded by Pharos, but it may still be useful with open-source or weaker models.

\subsubsection{Design}

\begin{figure}[ht]
\centering
\begin{tikzpicture}[x=0.985cm,y=1cm,
    role/.style={box,text width=36mm,minimum height=13mm,inner sep=1.5mm,line width=0.8pt},
    wk/.style={box,draw=good!70!black,fill=good!14,text width=20mm,minimum height=24mm,
      inner sep=1mm,line width=0.8pt,anchor=north},
    lm/.style={draw=black!45,dashed,fill=white,rounded corners=1mm,font=\footnotesize,
      inner xsep=0.5mm,inner ysep=1mm},
    elab/.style={lab,font=\footnotesize\itshape,align=center}]
  \node[role,draw=black!45,fill=black!6] (math) at (2.0,0)
    {\textbf{Mathematician}\\{\footnotesize\color{black!70} poses a problem in natural language}};
  \node[role,pbox,fill=tint,line width=0.8pt] (main) at (7.6,0)
    {\textbf{Main agent}\\{\footnotesize\color{black!70} global planning and orchestration}};
  \node[role,draw=yellow!55!black,fill=yellow!14] (pro) at (14.05,0)
    {\textbf{GPT-5.5 Pro}\\{\footnotesize\color{black!70} high-level strategy}};
  \draw[flow,{Stealth[length=2mm]}-{Stealth[length=2mm]}] (math) -- node[above,elab]{talk} (main);
  \draw[flow,dashed,{Stealth[length=2mm]}-{Stealth[length=2mm]},color=yellow!55!black]
    (main) -- node[above,elab]{consult} (pro);
  \draw[draw=good!70!black,dashed,fill=good!5,rounded corners=2mm,line width=0.7pt]
    (0,-1.75) rectangle (9.9,-5.5);
  \node[anchor=north west,font=\small,color=good!60!black] at (0.1,-1.8)
    {\textbf{Worker swarm}\; {\footnotesize modified Codex agents $\times N$}};
  \foreach \x/\eff [count=\i] in {1.32/high,3.74/high,6.16/xhigh,8.58/xhigh} {
    \node[wk] (w\i) at (\x,-2.65) {\textbf{worker}\\{\footnotesize\color{black!70} \eff\ effort}\\[9mm]\mbox{}};
    \node[lm] at ($(w\i.south)+(0,5.5mm)$) {local memory};
  }
  \node[role,abox,fill=warmtint,line width=1pt,minimum height=26mm,text width=30mm] (ver) at (14.35,-3.6)
    {\textbf{Verifier}\\{\footnotesize\color{black!70} stateless service}\\{\footnotesize\color{black!70} sole authority on correctness}};
  \draw[flow,color=good!70!black] (9.9,-3.0) -- node[above,elab]{submit statement\\+ proof} (ver.west|-0,-3.0);
  \draw[aflow,dashed] (ver.west|-0,-4.3) -- node[below,elab]{feedback} (9.9,-4.3);
  \draw[pflow] (main.south) -- node[right,elab]{assign $\cdot$ start $\cdot$ stop} (main.south|-0,-1.75);
  \draw[draw=black!45,dashed,fill=black!3,rounded corners=2mm,line width=0.7pt]
    (0,-6.5) rectangle (16.0,-9.0);
  \node[anchor=north west,font=\small\bfseries,color=black!60] at (0.1,-6.55) {Shared storage};
  \node[box,text width=66mm,minimum height=12mm,line width=0.8pt] (gm) at (3.9,-8.05)
    {\textbf{global memory}\\{\footnotesize\color{black!70} plans $\cdot$ findings $\cdot$ dead ends}};
  \node[box,draw=primary,fill=primary,text=white,text width=66mm,minimum height=12mm] (fg) at (12.1,-8.05)
    {\textbf{fact graph}\\{\footnotesize verified facts $\cdot$ the single source of truth}};
  \draw[flow,color=good!70!black,{Stealth[length=2mm]}-{Stealth[length=2mm]}]
    (3.9,-5.5) -- node[right,elab]{read $\cdot$ write} (3.9,-6.5);
  \draw[aflow] (ver.south) -- node[right,elab,pos=0.3]{writes\\verified fact} (ver.south|-fg.north);
\end{tikzpicture}
\caption{The Danus architecture (redrawn from \cite{liu2026danus}).}
\label{fig:danus}
\end{figure}
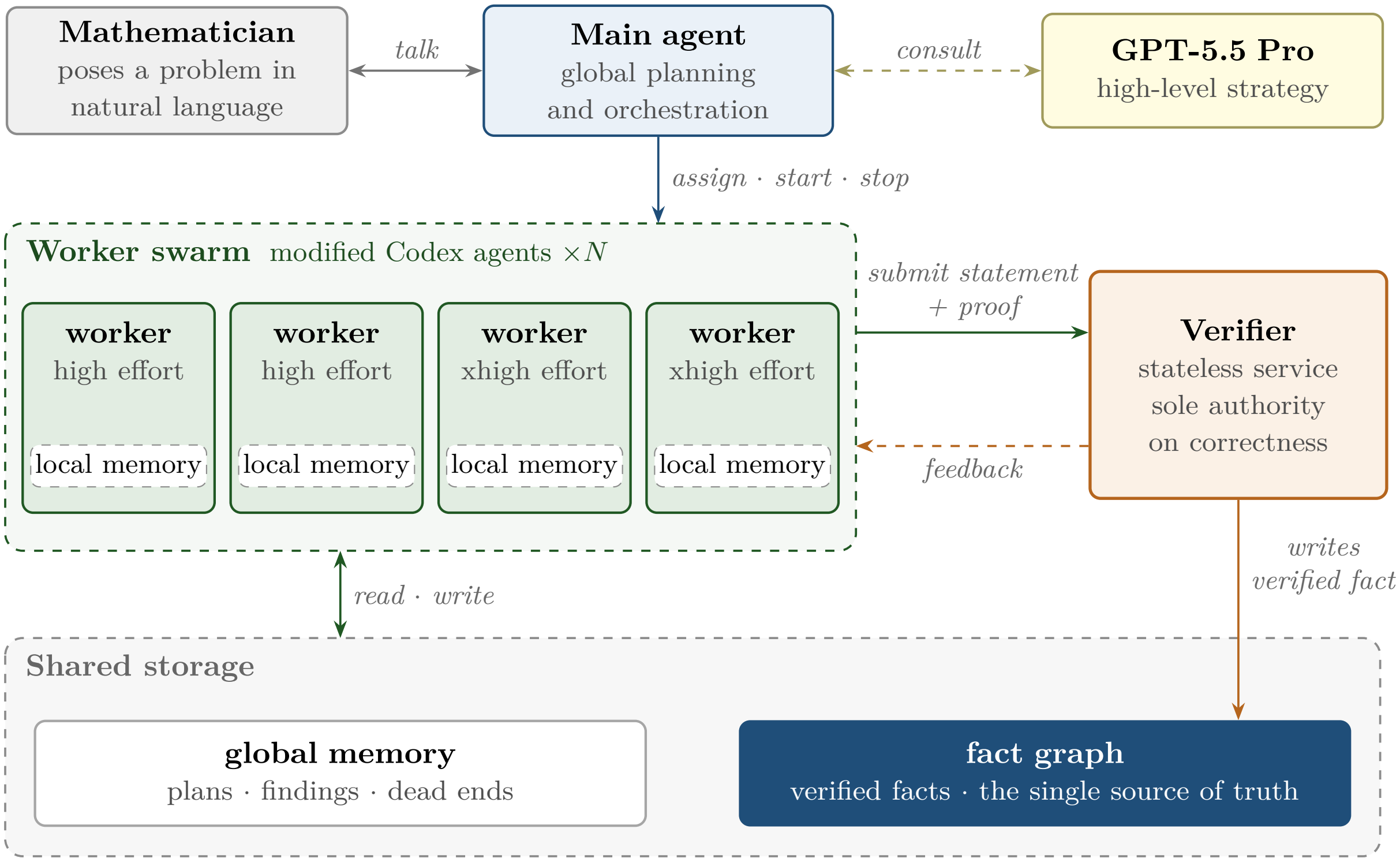
\paragraph{Agentic workflow}
Danus mimics a research group (\cref{fig:danus}).
\begin{itemize}
  \item The \textbf{main agent} is the group leader (Claude Code with Opus~4.8 in the paper's experiments). It plans routes, assigns and monitors work, and reports to the human. It does no mathematics.
  \item A \textbf{senior expert} (GPT-5.5 Pro) is consulted at most about once an hour on a distilled summary of the state.
  \item The \textbf{verifier} is stateless: a fresh instance judges each submitted claim.
  \item \textbf{Workers} are Rethlas generation agents (Codex with GPT-5.5), three to nine per project. Each owns one small problem at a time, assigned by the main agent.
\end{itemize}
\paragraph{The fact graph} The deliverable is no longer a single document but a directed acyclic graph, whose nodes are \emph{facts} (a statement with a proof accepted by the verifier) and whose edges are logical dependencies.

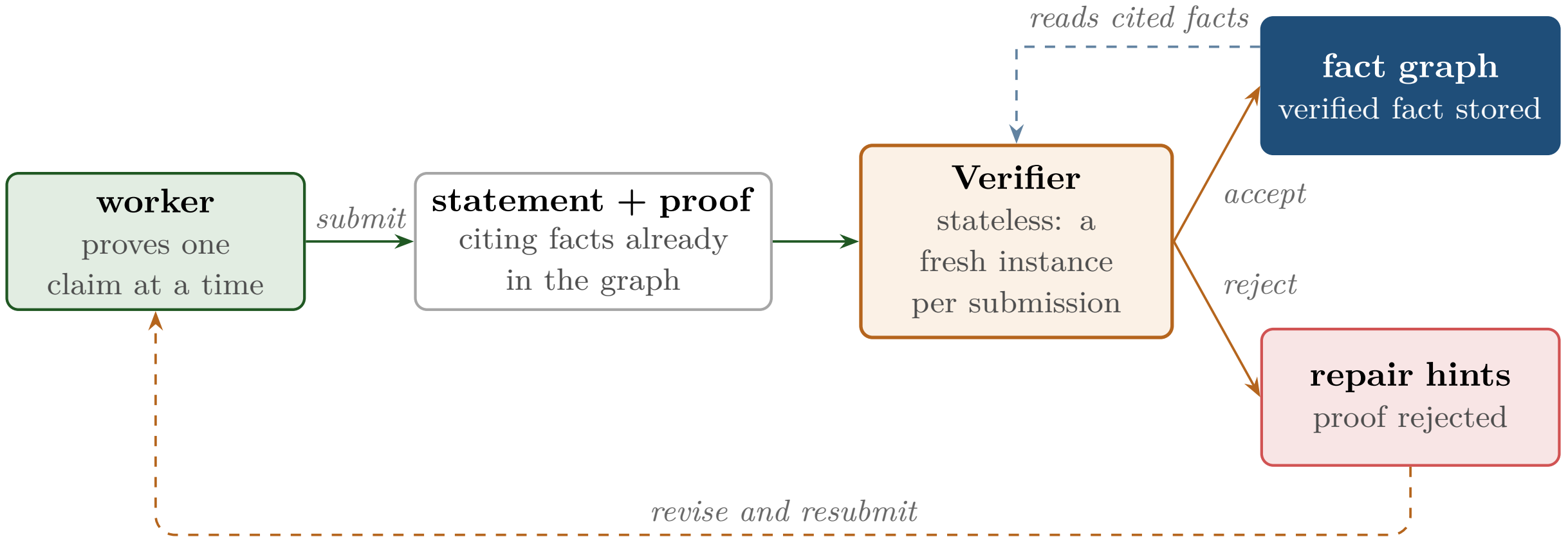
\begin{figure}[ht]
\centering
\begin{tikzpicture}[x=1cm,y=1cm,
    nd/.style={box,text width=28mm,minimum height=14mm,inner sep=1.3mm,line width=0.8pt},
    elab/.style={lab,font=\footnotesize\itshape,align=center}]
  \node[nd,draw=good!70!black,fill=good!14] (w) at (1.5,0)
    {\textbf{worker}\\{\footnotesize\color{black!70} proves one claim at a time}};
  \node[nd,text width=34mm] (sp) at (6.0,0)
    {\textbf{statement + proof}\\{\footnotesize\color{black!70} citing facts already\\in the graph}};
  \node[nd,abox,fill=warmtint,line width=1pt] (v) at (10.35,0)
    {\textbf{Verifier}\\{\footnotesize\color{black!70} stateless: a fresh instance per submission}};
  \node[nd,draw=primary,fill=primary,text=white] (fg) at (14.4,1.6)
    {\textbf{fact graph}\\{\footnotesize verified fact stored}};
  \node[nd,draw=bad!80,fill=bad!12] (rh) at (14.4,-1.6)
    {\textbf{repair hints}\\{\footnotesize\color{black!70} proof rejected}};
  \draw[flow,color=good!70!black] (w) -- node[above,elab]{submit} (sp);
  \draw[flow,color=good!70!black] (sp) -- (v);
  \draw[aflow] (v.east) -- node[below right=-1mm and 0mm,elab,pos=0.4]{accept} (fg.west);
  \draw[aflow] (v.east) -- node[above right=-1mm and 0mm,elab,pos=0.4]{reject} (rh.west);
  \draw[flow,dashed,color=primary!70] (fg.west|-0,2.0) -| node[pos=0.25,above,elab]{reads cited facts} (v.north);
  \draw[aflow,dashed,rounded corners=2mm] (rh.south) -- ++(0,-0.7) -| node[pos=0.25,above,elab]{revise and resubmit} (w.south);
\end{tikzpicture}
\caption{A Rethlas loop inside Danus.}
\label{fig:verifyloop}
\end{figure}

\paragraph{Why a fact graph?}
The fact graph serves several purposes:
\begin{enumerate}
  \item it \textit{decomposes} the problem into small steps, so that each worker faces an easier task;
  \item it is a \textit{communication channel} of verified results between workers;
  \item it gives \textit{cleaner context}: a worker loads only the facts it needs, by statement, and reads a proof only when its technique is needed;
  \item mistakes are easier to fix: \textit{revoking} a fact removes it and every fact that depends on it;
  \item it \textit{scales} compute on generation: each worker, together with the stateless verifier, forms a Rethlas loop (\cref{fig:verifyloop}), and since the deliverable is modular, parallel edits do not collide.
\end{enumerate}

\paragraph{Three-tier memory} Danus separates truth from awareness:
\begin{center}\small
\begin{tabularx}{\textwidth}{@{}l>{\raggedright\arraybackslash}X>{\raggedright\arraybackslash}p{30mm}>{\raggedright\arraybackslash}p{30mm}@{}}
\toprule
\textbf{Tier} & \textbf{Contents} & \textbf{Status} & \textbf{Human analogue}\\
\midrule
Fact graph & verified facts and their dependencies & the only truth & the group's verified draft (e.g. an Overleaf Project)\\
Global memory & plans, examples, counterexamples, dead ends & shared, \emph{not} truth & ideas shared in the group\\
Local memory & each worker's notes and logs & scratch paper & personal notes\\
\bottomrule
\end{tabularx}
\end{center}

\paragraph{Compaction at every boundary} Raw worker logs are compacted into global memory, global memory into the main agent's checkpoint
summaries, and those into reports for the human. A worker sees only its task, the compacted global memory, and the facts it retrieves from the fact graph. What any single agent reads stays bounded.

\paragraph{Separation of powers, enforced by tools} At scale, separation of
roles becomes separation of powers. In Danus it is enforced by which MCP tools each role can call, rather than by instructions in Markdown files:
\begin{center}\small
\begin{tabular}{@{}lccc@{}}
\toprule
& \textbf{main} & \textbf{worker} & \textbf{verifier} \\
\midrule
read facts & \yes & \yes & \yes \\
write memory & \yes & \yes & \no \\
submit a fact (\texttt{fact\_submit}) & \no & \yes & \no \\
decide truth & \no & \no & \yes \\
\bottomrule
\end{tabular}
\end{center}
In particular, the main agent cannot submit facts to the fact graph.

\paragraph{Scaling generation in both width and depth} \emph{Width}: the fact graph splits the proof into units a worker can own, so contributions accumulate, and the main agent diversifies directions. \emph{Depth}: per-fact verification and the memory hierarchy extend the reasoning horizon beyond one context window. 

\paragraph{Completion and paper writing}
Danus declares success when the original statement of the problem enters the fact graph. The main agent then writes a paper from the fact graph and submits it to the verifier; this write--verify loop continues until the verifier accepts the paper.

\subsubsection{Case study: tangent classes of matroid building sets}
The paper \cite{liu2026danus} contains six case studies; here we discuss one of them, tangent classes of matroid building sets, in more detail.
\paragraph{Problem statement}
Let $M$ be a loopless matroid (not necessarily realizable) of rank $d+1$ and
$\mathcal G$ a building set. Construct a class $T_{M,\mathcal G}$ in the
combinatorial $K$-ring $K(M,\mathcal G)$ such that (i) if $M$ is realizable,
$T_{M,\mathcal G}$ is the tangent class of the wonderful compactification;
(ii) a Hirzebruch--Riemann--Roch formula
$\dim A^i(M,\mathcal G) = (-1)^i \chi(\wedge^i T^\vee)$ holds; and (iii)
the Chern numbers satisfy $\deg(c_k(T)\alpha^{d-k}) \ge \binom{d+1}{k}$.

This is the main result of Cheng's paper \cite{cheng2026tangentclassesmatroidbuilding}, obtained independently.
\paragraph{Comparison with base models}
With the same verbatim prompt, GPT-5.5 Pro (web) produced no progress and no strategy plan; three runs of Rethlas produced proof drafts, none of which passed the internal verification agent (0/3); Danus produced a verified solution \cite{cheng2026tangent}.
\paragraph{How we selected this problem}
By coincidence. Over lunch on California Avenue in May 2026, Ronnie Cheng told me that he had just completed his paper \cite{cheng2026tangentclassesmatroidbuilding}, and that GPT-5.5 Pro could not prove its main results. He generously agreed to let us try our systems on them, and delayed posting his paper to arXiv. GPT-5.5 Pro and Rethlas both failed; Danus re-proved the main results \cite{cheng2026tangent}.

\begin{figure}[ht]
\centering
\includegraphics[width=0.32\linewidth]{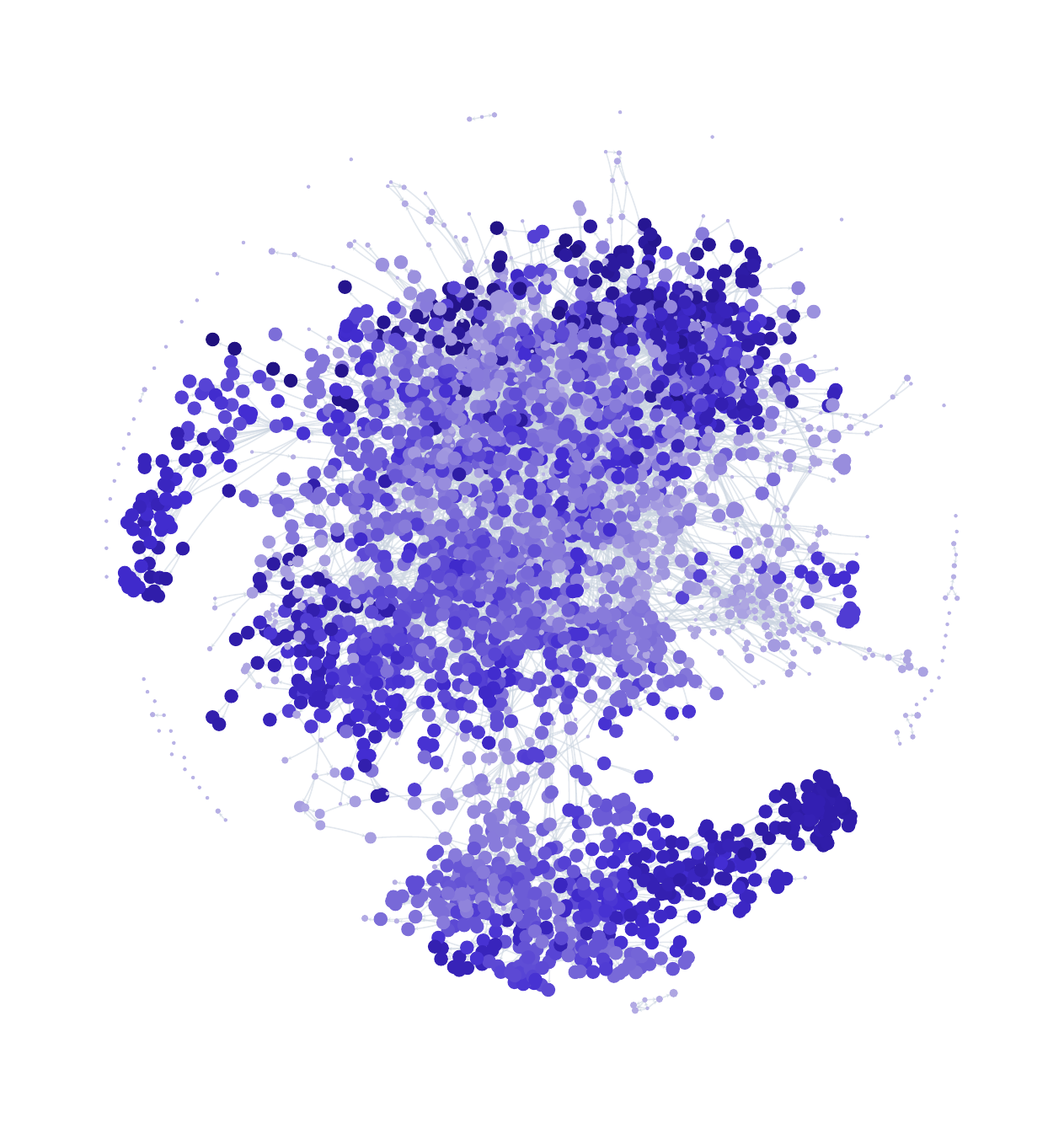}
\caption{The fact graph of the tangent-class run \cite{liu2026danus}.}
\label{fig:factgraph}
\end{figure}

\paragraph{Experiment setup}
The run used 7 workers (3 at \texttt{xhigh} and 4 at \texttt{high} effort). It took about five days, with about 50 active hours and 13 GPT-5.5 Pro consultations.
\paragraph{Danus scales width and depth of generation}
The fact graph (\cref{fig:factgraph}) has 3{,}157 verified facts and 8{,}616 dependency edges, with chains up to 54 facts deep. Only 499 facts are needed for the main theorem, so about six times more facts were proved than the final proof uses.

\paragraph{Compute can buy intelligence} The logs of the Rethlas runs show that Rethlas conflated nearby notions (realizable
matroids, wonderful models, toric varieties; integral versus rational
$K$-theory). This suggests that mathematical subtlety, such as the difference between nearby notions, can be overcome by scaling compute in width and depth: thinking longer and more carefully.

\subsubsection{Limitations}
Danus has several limitations.
\begin{itemize}
    \item \textit{Occasional verifier false positives.} Once compute on generation is scaled, compute on verification must be scaled too, so that errors in the fact graph do not ruin the run. This motivates OSVerify (\cref{sec:verification}).
    \item \textit{Outdated by stronger base models.} When the base model improves, each worker should  be assigned harder tasks rather than tiny steps. This motivates Pharos (\cref{sec:pharos}).
    \item \textit{Writing and communication with human experts.} The papers and reports are hard for humans to read quickly, because of over-compressed steps, unconventional notation, and so on.
\end{itemize}

\subsection{Pharos: Danus for stronger models}\label{sec:pharos}
Pharos is fully open source, with code at \url{https://github.com/frenzymath/Pharos}, a technical report at \url{https://frenzymath.com/news/pharos-preview/}, and a paper in preparation \cite{pharos}; see also \cite{gao2026harnesstalk}.

\subsubsection{Design}
\paragraph{Overview of comparison with Danus} Pharos inherits Danus's fact graph, workers, and verifier, but rebuilds the rest for GPT-6 Astra and other stronger models. Some components that
compensated for weaker models were removed. \Cref{tab:danus-pharos} lists selected changes.
\begin{table}[ht]
\centering\small
\renewcommand{\arraystretch}{1.2}
\begin{tabularx}{\textwidth}{@{}l>{\raggedright\arraybackslash}p{45mm}>{\raggedright\arraybackslash}X@{}}
\toprule
& \textbf{Danus} (July 2026) & \textbf{Pharos v3.1} (Sept.\ 2026) \\
\midrule
Main agent & Claude Code, Opus 4.8 &  Codex, \texttt{ultra} effort \\
Consultant & GPT-5.5 Pro consulted every 1--2\,h & consultant removed \\
Planning & periodic summaries & \emph{route board}\\

Writing & one write-verify loop & outline $\to$ per-section writer--verifier pairs $\to$ whole-paper verifier; distilled expert writing style \\
Computation & no large numerical computations & allowed, with compute caps and a progress monitor\\
Workers & per-round rhythm & multi-round work on one result \\
Memory & full scans of memory & incremental index for memory search \\
Human input & ad hoc & recorded verbatim \\
Models & GPT-5.5 workers & GPT-6 Astra workers \\
\bottomrule
\end{tabularx}
\caption{From Danus to Pharos.}
\label{tab:danus-pharos}
\end{table}

\paragraph{The route board}

\begin{figure}[ht]
\centering
\begin{tikzpicture}[x=1cm,y=1cm,font=\small,
    rstep/.style={box,font=\small\bfseries,minimum width=4.5cm,minimum height=1.0cm}]
  \node[fill=primary,text=white,rounded corners=4pt,minimum width=3cm,
    minimum height=1.1cm,font=\normalsize\bfseries] (B) at (0,0) {Route board};
  \node[rstep] (s1) at (-5.2,1.2) {1.\ Survey known methods};
  \node[rstep] (s2) at (0,2.4) {2.\ Propose routes};
  \node[rstep] (s3) at (5.2,1.2) {3.\ Test special cases};
  \node[rstep,fill=tint,draw=primary] (s4) at (5.2,-1.2) {4.\ Work out details};
  \node[rstep,fill=warmtint,draw=accent] (s5) at (0,-2.4) {5.\ Judge progress};
  \node[rstep] (s6) at (-5.2,-1.2) {6.\ Find the obstruction};
  \foreach \s in {s1,s2,s3,s4,s5,s6} \draw[black!15,line width=.6pt] (B) -- (\s);
  \draw[pflow] (s1) -- (s2);
  \draw[pflow] (s2) -- (s3);
  \draw[pflow] (s3) -- (s4);
  \draw[pflow] (s4) -- (s5);
  \draw[pflow] (s5.west) -- (s6.south east);
  \node[anchor=north east,lab] at ($(s5.west)+(-.05,-.2)$) {continue $\cdot$ reassign $\cdot$ stop};
  \draw[pflow,dashed] (s6) -- node[left,lab,align=right]{refine into\\sub-routes} (s1);
  \node[pbox,fill=tint,font=\footnotesize] (fg) at (5.2,-2.9) {\textbf{Fact graph}: verified, reusable};
  \draw[flow,color=good] (s4) -- (fg);
\end{tikzpicture}
\caption{How the Pharos route board is updated
\cite{gao2026harnesstalk}.}
\label{fig:routeboard}
\end{figure}
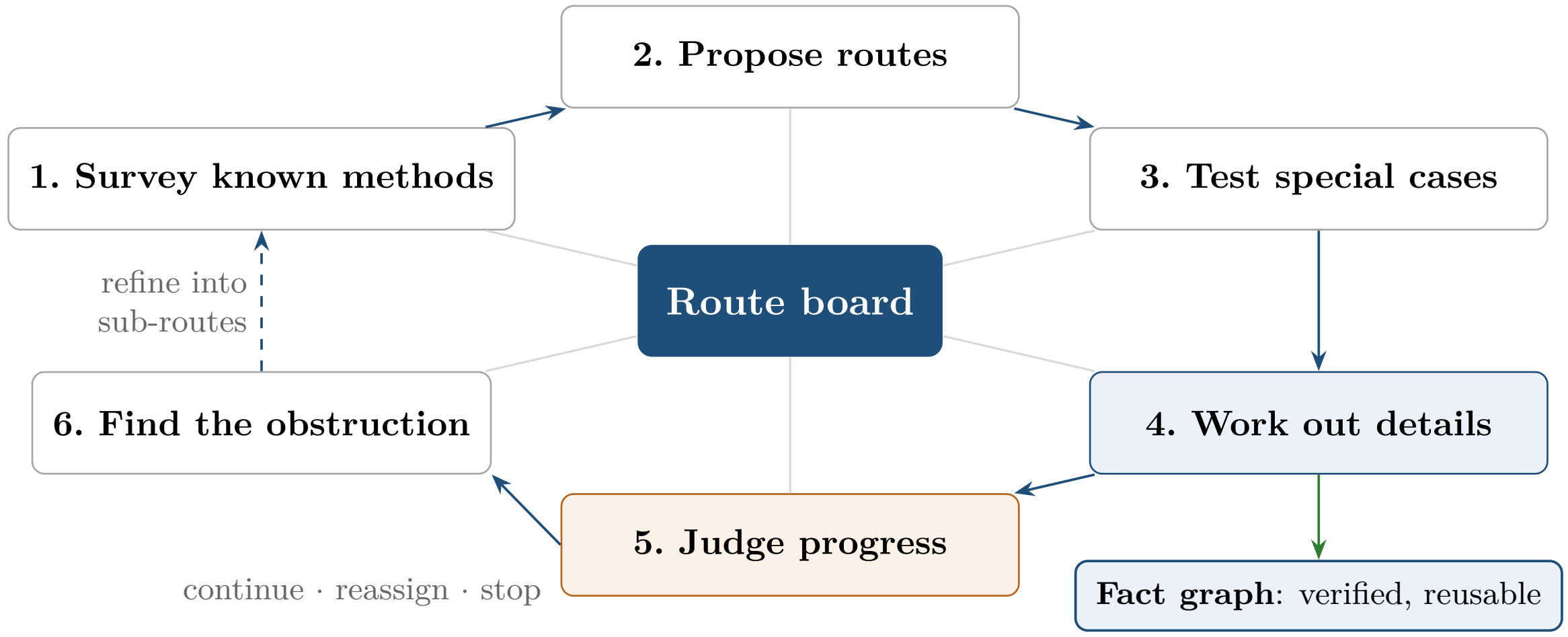

The main new structure is a \emph{route board} (the file \file{ROUTES.md}),
the main agent's live record of strategies and plans.
\Cref{fig:routeboard} shows the update cycle.

\begin{description}
  \item[Known routes first.] Every technique from the literature is
  tried before new routes are invented.
  \item[Live status.] The board states the current status and what to try next. This addresses the long-context confusion of \cref{sec:what-goes-wrong}: each worker sees only its own route's
  briefing.
\end{description}

\paragraph{Other improvements}

\begin{itemize}
  \item \textbf{Deeper facts.} Danus workers spent on the order of ten minutes
  on each small fact; Pharos workers spend one to two hours on harder facts, which suits GPT-6 Astra.
  \item \textbf{Helper agents} digest literature, scout routes, and assess progress. 
  \item \textbf{Computation.} Workers run code and numerical experiments with budget caps.
  \item \textbf{Working with experts.} An expert can follow a run and steer it through the main agent.
  \item \textbf{Efficiency.} Lossless context compression lowers cost and raises throughput.
  \item \textbf{Improved writing pipeline.} Once the target is a fact, the main agent writes an outline and assigns sections; each
  section writer revises against its own verifier; a whole-paper verifier checks the assembled argument; the paper is then polished and built as a PDF.
\end{itemize}
\paragraph{New role: operations agent}
The implementation adds an \emph{operations agent} at the repository root, which deploys and monitors several projects and implements user customizations.

\subsubsection{Case study 1: an algebraic proof of the Hartshorne conjecture}
The first case study is an algebraic, characteristic-zero proof of the Miyaoka--Mori criterion, which reproves the Hartshorne conjecture \cite{dong2026jets}. How Pharos was used is described in \cite[Appendix~B]{dong2026jets}.
\paragraph{Problem statement}
For a smooth projective variety $X$ in characteristic zero and a curve
$f\colon C \to X$, does $-K_X \cdot f_*[C] > 0$ alone give a rational curve
through every point of $f(C)$, under a strict degree bound?
\paragraph{Comparison with base model}
Du, Guo, and Xie tried GPT-6 Astra Pro in the web interface and did not obtain a solution, even after many rounds of mathematical guidance. In contrast, Pharos with GPT-6 Astra solved the problem fully automatically, in 19\,h\,57\,min with no human intervention \cite{dong2026jets}.
\paragraph{Experiment setup}
The main agent ran at \texttt{ultra} effort, 3 workers at \texttt{xhigh} and 4 at
\texttt{max}. The total cost was under twice the weekly Pro usage limit.
\paragraph{Deeper facts and exploration}
The fact graph has 273 verified facts (one revoked and re-proved), of which only 4 are in the final proof.
\paragraph{Different routes and interactions}
There were 7 routes, refined into 30 sub-routes, and facts were reused across routes. The final proof combines routes R3 and R6 through a new sub-route, R3-6.
\paragraph{Human hints can save compute}
With expert hints at launch, Pharos found the same proof in 11\,h\,38\,min.
\subsubsection{Case study 2: creating a framework to unify two problems}
This case study was carried out by A.~Shelley and H.~Bosch, using Pharos
\cite{shelley2026collision}. The \textit{loc.cit} cites Danus because at that time Pharos was not named.

\paragraph{Background} Shelley had several related lower bounds for two
problems: stationary geodesic networks, and Newtonian point billiards. Bosch
had proved a superexponential bound in his work on stationary geodesic
networks \cite{bosch2026regularity}. They wanted to place the two problems
in a common light, and to optimize the constructions.

\paragraph{From one bound to a framework} On its first attempt, Pharos did
not improve the construction. Given Bosch's paper and the new lower bound, it
found, after about an hour, a proof that improved Bosch's bound from
superexponential to exponential. Bosch then suggested that the two problems
might fit a common framework. Given
this as a prompt, Pharos worked for between half an hour and an hour and
produced a manuscript unifying the two problems under a common definition,
using the formalism of Hamiltonian mechanics. The name \emph{collision systems} came from a human; the definition that makes the framework work came from the agent.

\paragraph{Optimizing the bounds} Asked next to optimize the bounds, Pharos
improved the lower bound on the optimal base from $2^{1/3}$ to $3^{1/4}$
(under a relaxation of the problem), and the upper bound from $8$ to $3.98$.
The authors then spent several days cleaning up the proof with Claude Fable,
reducing its complexity and improving its clarity, before posting it on
Hexagon \cite{shelley2026collision}.

\paragraph{Comparison with base model} According to the authors, Claude Fable~5.1 alone failed to find a common framework for the two problems.

\makeatletter
\renewcommand\subsubsection@cntformat{\thesubsubsection.}
\makeatother
\let\subsubsection\origsubsubsection
\let\paragraph\origparagraph
\setcounter{secnumdepth}{3}
\addtocontents{toc}{\protect\setcounter{tocdepth}{2}}

\section{Agents for verification: OSVerify}\label{sec:verification}\label{sec:osverify}
OSVerify is fully open source, with a paper \cite{nie2026osverify}, a website
at \url{https://osverify.github.io}, and code at \url{https://github.com/OSVerify}.

\subsection{Motivation and objective}

\subsubsection{For human experts}
Verification has become the expensive step. Papers can now be generated much faster than they can be refereed, so we need to reduce the human effort spent on refereeing.

\subsubsection{For reasoning agents}
Once an agentic system scales up generation (as Danus and Pharos do), a weak
verifier ruins the run: false lemmas pile up in the fact graph like a
denial-of-service attack. Test-time compute should therefore be scaled on the
verifier as well.

\subsubsection{Objective: a stricter reviewer}
We want a stricter verifier that reliably rejects wrong solutions while rarely
rejecting correct ones. Formally, a verifier $V(x,y) \in \{\textsc{Reject},
\textsc{Pass}\}$ judges a solution $y$ to a problem $x$, and we want it to be
one-sided in the following sense.
\begin{definition}[One-sided verification]
Maximize the true negative rate $\TNR$ subject to $\FNR \le \tau$ for a small
tolerance $\tau$, where
\[
\TNR = \Pr[\textsc{Reject} \mid y \text{ invalid}], \qquad
\FNR = \Pr[\textsc{Reject} \mid y \text{ valid}].
\]
\end{definition}

\subsection{Design: a propose--validate loop}

\begin{figure}[ht]
\centering
\makebox[\linewidth][c]{\includegraphics[width=1.12\linewidth]{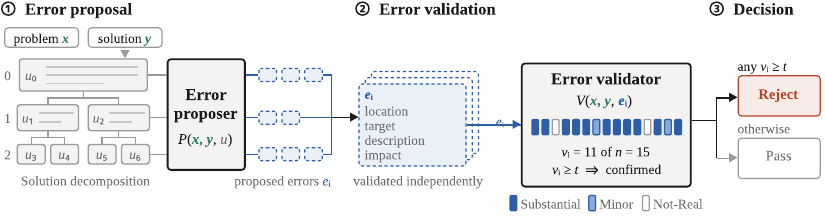}}
\caption{OSVerify \cite{nie2026osverify}: an error proposer decomposes the
solution and proposes errors; an error validator judges each proposed error.}
\label{fig:osv}
\end{figure}

\subsubsection{Agentic workflow}
OSVerify consists of two agents (\cref{fig:osv}).
\begin{description}[font=\normalfont\itshape]
  \item[Error proposer (goal: coverage).] Decomposes the solution
  hierarchically and proposes candidate errors.
  \item[Error validator (goal: no false alarms).] Classifies each candidate
  error as \textsc{Not-Real}, \textsc{Minor}, or \textsc{Substantial}.
\end{description}

\subsubsection{Reasons for the design}
Several observations motivate this design.
\begin{itemize}
  \item \textit{Mimic the human workflow.} A human referee first
  \emph{suspects} a step and then \emph{checks} whether it is a real, fatal
  error. OSVerify assigns a separate agent to each task.
  \item \textit{Errors are easier to confirm than to find.} LLMs that fail to
  locate a reasoning error can often correct it once told where it is
  \cite{tyen2024llms}.
  \item \textit{Context engineering.} The proposer examines one piece at a
  time while the whole solution stays available; each validator call
  concentrates on a single error.
  \item \textit{Scale compute.} The generate--verify principle of Rethlas is
  applied \emph{inside} verification. The resulting second-order
  generate--verify loop scales compute and yields a stricter referee.
\end{itemize}

\subsection{Comparison with LLM-as-a-judge and majority voting}

\subsubsection{Benchmarks}
The benchmarks are Putnam100 (58 invalid and 50 valid solutions), Hard2Verify
(H2V) \cite{pandit2026hard2verify} (163 invalid, 37 valid), and
FirstProof-Eval, built from the twenty First Proof problems
\cite{abouzaid2026firstproof} (59 invalid, 35 valid). H2V and First Proof were
released after both backbones.

\subsubsection{Backbones}
The backbones are deliberately weak (\texttt{gpt-oss-120b} and
\texttt{gpt-5-mini}), to simulate problems beyond the model's reliable
ability.

\subsubsection{Training the validator}
The validator can also be trained by reinforcement learning (GRPO with LoRA)
on a dataset of labeled (problem, solution, error) triples; First Proof was
not used in this dataset. See \cite[Appendix~D]{nie2026osverify} for details.

\subsubsection{Results}
\Cref{tab:osv} compares OSVerify with LLM-as-a-judge and majority voting.

\begin{table}[ht]
\centering\small
\renewcommand{\arraystretch}{1.2}
\begin{tabular}{@{}lccc@{}}
\toprule
\textbf{Method} (\texttt{gpt-oss-120b}) & \textbf{Putnam100} & \textbf{H2V} & \textbf{FirstProof-Eval} \\
\midrule
LLM-as-a-judge & 0 & 54.4 & 0 \\
Voting (16 samples) & 38.4 & 69.2 & 35.9 \\
OSVerify & 70.4 & \textbf{84.7} & 58.0 \\
\quad + RL-trained validator & \textbf{89.6} & 83.2 & \textbf{62.3} \\
\bottomrule
\end{tabular}
\caption{TNR (\%) at $\FNR \le 5\%$ \cite{nie2026osverify}. A $0$ means that
only the always-pass rule meets the constraint.}
\label{tab:osv}
\end{table}

\subsubsection{Scaling the proposer and the validator}
\begin{figure}[ht]
\centering
\makebox[\linewidth][c]{\includegraphics[width=1.12\linewidth]{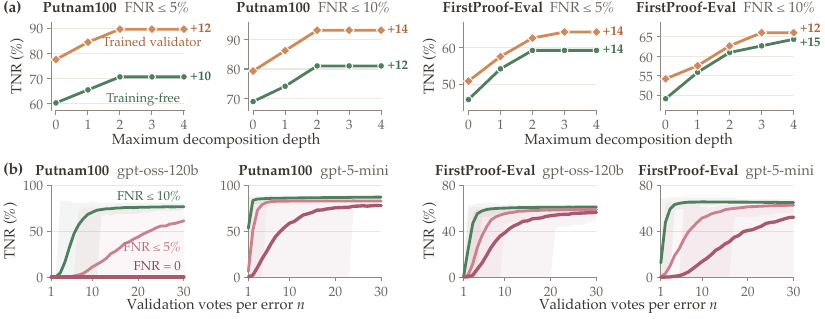}}
\caption{Scaling compute on verification \cite{nie2026osverify}: (a) proposer
depth; (b) validator votes.}
\label{fig:osvscaling}
\end{figure}

Both components benefit from more compute (\cref{fig:osvscaling}). Deeper
decomposition raises TNR, with most of the gain reached by depth two. More
validator votes also raise TNR, and stricter FNR constraints need more votes.
At every matched token budget OSVerify beats voting, and the RL-trained
validator reaches the same TNR with about a quarter of the compute.

\subsection{Comparison with frontier models: case studies in research mathematics}

With frontier backbones, OSVerify has been used to referee research papers.

\subsubsection{AI-generated papers}
OpenAI released eleven preprints, each claiming progress on a central problem
(Kakeya in dimension~4, planar unit distances, the trilinear Hilbert
transform, quantum geometric Langlands, \ldots). OSVerify found a nontrivial
mistake in each of them. The mistakes had been overlooked by OpenAI's internal
review and were confirmed by experts, although the main conclusions survived
revision. See \url{https://osverify.github.io/case-studies.html}.

\subsubsection{Kamnitzer's ICM 2022 conjecture}
Q.~Tang wrote an AI-assisted manuscript proving that the Mirković--Vilonen
basis and the theta basis coincide \cite{tang2026mvtheta}. In round~1, GPT-6
Pro drafted the proof and Codex reviewed it; OSVerify caught a major error. In
round~2, GPT-6 Pro repaired the proof, and three further reviews by ChatGPT and
Codex found nothing; OSVerify caught another error that broke the proof.

\subsubsection{The Manin constant conjecture}
Z.~Chi, S.~Liu and P.~Nie obtained an AI-assisted proof of the Manin constant
conjecture \cite{chi2026manin}, which arose in the study of the Birch and
Swinnerton-Dyer conjecture and had been open for over 55 years. Twelve
independent rounds of AI checks (Claude Opus~5.5, Claude Fable~5.1, and Codex
with GPT-6 Astra) found no error in an early draft. OSVerify (with GPT-5.6
Sol) found one and suggested fixes for the three places it affected. See
\cite[Appendix~A]{chi2026manin}.

\subsection{OSVerify as a meta-algorithm}

\subsubsection{For verification}
Any verifier that raises errors, such as multi-agent debate
\cite{du2024multiagentdebate} or pseudo-formalization
\cite{barkallah2026pseudoformalization}, can serve as the proposer, with
OSVerify's validator on top. On H2V at $0\%$ FNR, this raised the TNR of
multi-agent debate from $71.2\%$ to $80.8\%$, and that of pseudo-formalization
from $0$ to $92.6\%$.

\subsubsection{For problem solving}
OSVerify can replace the simple verifier in problem-solving agents such as
Rethlas, Danus, and Pharos; Rethlas-OSV (\cref{sec:use}) is an example.

\section{How to use and contribute as a mathematician}\label{sec:use}
This section first explains how to use the agents introduced in this paper,
and then how to contribute to improving them.
It focuses on the systems introduced in this paper. For general advice on
using AI in mathematical research, we refer to other sources, for example
Etingof's guide for young mathematicians \cite{etingof2026aiuse}.

The agents are designed for Linux and macOS, and the instructions below are
written for these two systems. On Windows one can use WSL or a virtual
machine, and ask a coding agent for help, but performance is not guaranteed.

\subsection{Choosing a tool: don't use a sledgehammer to crack a nut}\label{subsec:choose-tool}

\subsubsection{Empirical suggestions}
\Cref{tab:use} gives the author's empirical estimates as of November 2026.
Runtimes and quotas vary widely with the problem and the model.

\begin{table}[ht]
\centering\small
\renewcommand{\arraystretch}{1.25}
\begin{tabularx}{\textwidth}{@{}l>{\raggedright\arraybackslash}X>{\raggedright\arraybackslash}p{36mm}l@{}}
\toprule
\textbf{System} & \textbf{Use it for} & \textbf{Typical runtime} \\
\midrule
Chat (GPT Pro, web) & discussion, small questions, literature & minutes \\
Codex / Claude Code (\texttt{/goal}) & a problem you expect to be solvable & 30 min -- 3 h \\
Rethlas & a problem whose solution may have a simple logical structure & 1--5 h  \\
Danus / Pharos & a complicated problem needing many lemmas & about a weekly Pro quota per 12 h  \\
OSVerify & checking a paper or a proof & 1--5 h ($<$10\% of a weekly quota) \\
Rethlas + OSVerify & solving with a stricter verifier & 1--5 h ($<$20\% of a weekly quota with 1 branch) \\
\bottomrule
\end{tabularx}
\caption{Which tool for which job.}
\label{tab:use}
\end{table}

\subsubsection{Danus or Pharos}
With stronger models, such as GPT-6 Astra, GPT-6.1 Sol, or Claude Opus~5.5,
Pharos is strongly preferred: Danus has been observed to limit the
performance of stronger models. With weaker models, such as GPT-5.5 or
open-source models, Danus may be the better choice, and a strong model can
serve as its external consultant agent.

\subsubsection{OSVerify as the verifier}
Pharos can also be combined with OSVerify, but this may be very expensive.

\subsection{Installing and running the agents via a coding agent}
The easiest way is to ask Codex or Claude Code to download, install, and run
these systems for you from source. Here are some template prompts.

\subsubsection{Install Rethlas}
\begin{promptbox}
Help me download Rethlas from \url{https://github.com/frenzymath/Rethlas}.
Check that it runs correctly on my laptop and that its MCP tools work.
\end{promptbox}

\subsubsection{Run Rethlas}
\begin{promptbox}
Help me run Rethlas, with GPT-5.6 Sol at max effort as the generation agent
and GPT-5.6 Sol at xhigh effort as the verification agent, on the following
problem: prove or disprove the BSD conjecture.
\end{promptbox}

\subsubsection{Pharos on macOS}
By default, Danus and Pharos run only on Linux. You can run them on a server
(for example an AWS, or FarmShare at Stanford; most servers run
Linux), or ask your coding agent to adapt them to macOS:
\begin{promptbox}
Help me install Pharos from \url{https://github.com/frenzymath/Pharos} and
adapt it to macOS. Check that it runs correctly on my laptop and that its MCP
tools work.
\end{promptbox}

\subsubsection{Desktop app: Rethlas-OSV}
Rethlas-OSV is a desktop app for macOS and Linux. It runs Rethlas with
OSVerify for problem solving, and OSVerify alone for refereeing papers.
Download it from \url{https://github.com/OSVerify/Rethlas-OSV/releases}; a
step-by-step guide is at \url{https://srliu.com/blog/rethlas-osv-guide/}. For
a first attempt we suggest one to three parallel branches and 20 rounds; you
can later resume the run with more rounds.

\subsection{Running and customizing the agents manually}
More detailed guidance is kept at \url{https://srliu.com/blog}.

\subsubsection{Prerequisites}
You need the Codex CLI (see
\url{https://learn.chatgpt.com/docs/codex/cli#getting-started}) and some
familiarity with the terminal (see \url{https://srliu.com/blog/terminal-guide/}).

\subsubsection{Rethlas}
The verifier runs as a local web service, and the generator is a Codex agent
that calls it. Using tmux is convenient but not required.

Installation:
\begin{shell}
git clone https://github.com/frenzymath/Rethlas.git && cd Rethlas

# 1. set up the verification agent
cd agents/verification && python3 -m venv .venv && source .venv/bin/activate
pip install -r requirements.txt && deactivate && cd ../..

# 2. set up the generation agent
cd agents/generation && python3 -m venv .venv && source .venv/bin/activate
pip install -r mcp/requirements.txt
\end{shell}
Usage:
\begin{shell}
# 1. start the verification service
cd Rethlas/agents/verification && source .venv/bin/activate
uvicorn api.server:app --host 0.0.0.0 --port 8091

# 2. in another shell: run the generation agent on a problem
cd Rethlas/agents/generation && source .venv/bin/activate
PROBLEM_FILE=data/my_problem.md ./tests/run_example.sh
\end{shell}
Write each problem in a Markdown file and set \texttt{PROBLEM\_FILE} to its
path. The environment variables \texttt{MODEL}, \texttt{REASONING\_EFFORT},
and \texttt{VERIFY\_URL} select the model, the effort level, and the
verifier. References can be placed in \file{data/<problem>.refs/}; Markdown
and LaTeX files are preferred over PDF.

The outputs are the logs in \file{logs/<id>/}, the notebook in
\file{memory/<id>/}, and the blueprint in \file{results/<id>/blueprint.md},
which is renamed \file{blueprint\_verified.md} once it passes verification.

\subsubsection{Danus}
Clone the repository \url{https://github.com/frenzymath/Danus}, and then run
the following from its root.
\begin{shell}
# 1. install the toolchain (Node, Python venv, Codex CLI) into runtime/
bash scripts/bootstrap.sh

# 2. copy the config templates and fill in YOUR keys (never committed)
cp config/danus.env.example config/danus.env
cp config/codex.env.example config/codex.env   # OpenAI-compatible endpoint

# 3. health check, then start the verify service (required for proving)
bash scripts/doctor.sh
bash scripts/services.sh up verify

# 4. start Codex at the repository root; the first run calls `initialize`.
#    The flag below lets the main agent act without approval prompts, with
#    your shell privileges: run Danus on an isolated, disposable host, and
#    read docs/security-and-trust.md first.
codex --dangerously-bypass-approvals-and-sandbox
\end{shell}

\subsubsection{OSVerify}
Clone the repository \url{https://github.com/OSVerify/OSVerify}. You need
Python~3.12+, \texttt{uv}, an authenticated Codex CLI, and XeLaTeX or
LuaLaTeX for PDF reports.
\begin{shell}
uv sync
uv run osverifier demo --run-dir runs/demo     # offline demo, no model calls

# review a solution
uv run osverifier verify --question question.md --solution solution.md \
  --request-id problem-001 --model MODEL_ID --decompose \
  --run-dir runs/problem-001

# review a paper
uv run osverifier verify --review-type paper --decompose --question scope.md \
  --solution paper.tex --request-id paper-v1 --model MODEL_ID \
  --run-dir runs/paper-v1

# solve with Rethlas, using OSVerify as verifier
uv sync --extra rethlas
uv run osverifier-rethlas --problem problem.md --run-dir runs/my-problem \
  --generation-branches 4 --max-rounds 10

uv run osverifier-gui                          # local GUI
uv run osverifier export --run-dir runs/problem-001 --format latex \
  --output report.tex
\end{shell}

\subsection{Pharos via its operator agent}
Pharos is designed to be interactive and to be managed through a coding-agent
CLI.

\subsubsection{Requirements}
Pharos runs on Linux with Python~3.10+, Git, tmux, Chromium (for reports),
and \texttt{pdflatex} or \texttt{tectonic} (for papers). It can easily be
adapted to macOS.

\subsubsection{Start the operator agent}
\begin{shell}
git clone https://github.com/frenzymath/Pharos.git && cd Pharos
codex --yolo                         # from the repository root
\end{shell}
You are now talking to the operator agent. Ask it to set up Pharos and
confirm that it works on your machine. Then describe your problem and state
your requirements in natural language; the operator agent creates and
launches the project for you. You will be asked to review the problem
statement and the initial guidance before research begins. Afterwards, you
can send guidance either directly to the project's main agent or through the
operator agent.

To customize Pharos, describe what you want to the operator agent and let it
implement the changes.

\subsection{Good practice}

\subsubsection{State the problem precisely}
Avoid abuse of notation. If the statement may be false, say ``prove or
disprove''.

\subsubsection{Supply references}
You may need to download copyrighted books yourself. Where possible, provide
LaTeX or Markdown files rather than PDFs.

\subsubsection{Choose the right tool}
See \cref{subsec:choose-tool}.

\subsubsection{Be transparent}
Disclose the use of AI systems in your papers, say which parts were
AI-generated, and take responsibility for correctness.

\subsection{How to contribute}\label{sec:contribute}
Harness design is a place where mathematical expertise contributes directly.
Contributions are very welcome; please reach out to us.

\subsubsection{Find failure modes}
Run a system on a problem you understand, find how and why it fails, and
report it.

\subsubsection{Find new applications}
Look for new ways for AI to contribute meaningfully to mathematics: solving
problems in new areas, building theories, improving writing or understanding,
and so on. For example, \cite{bosch2026regularity} is the first instance of
using Pharos to develop a framework connecting two problems.

\subsubsection{Build datasets and benchmarks}
These are important both for evaluating agentic systems and for training
models.

\subsubsection{Distill workflows into skills and orchestration}
Human experience provides intuition for harness design. Improve the harness
based on your own research experience.

\section{Future directions and applications in other domains}\label{sec:future}
\subsection{Open directions}

\subsubsection{Simplifying and improving the harness as models improve}
As models get stronger, components that compensate for their weaknesses
become unnecessary (\cref{sec:agents}). Which components should be removed,
and what more ambitious harnesses do stronger models make possible?

\subsubsection{From verification to understanding and writing}
How can AI help humans understand papers, and help produce clearer ones? This
is the subject of our ongoing work.

\subsubsection{Human--AI collaboration}
We need more effective ways for experts to steer long runs (as in Danus and
Pharos), and better writing and communication skills for agents.

\subsubsection{Creativity and theory building}
Current systems work largely within the ``convex hull of known techniques''
(not a well-defined notion!). How can we make them more creative, and even
able to build new theories?

\subsubsection{Evaluation}
How should research-level agentic systems be benchmarked rigorously? The
evaluation of agentic systems is still largely empirical.

\subsubsection{Self-improvement}
Can a harness recursively improve itself?

\appendix
\section{Glossary}\label{sec:glossary}
\begin{description}[style=nextline]
  \item[Blueprint] Rethlas's proof document: a sequence of lemmas with
  statements and proofs, ending with the main theorem.
  \item[Compaction] Replacing old context by a summary when the context window
  fills up.
  \item[Context window] The text the model can attend to at once.
  \item[Fact, fact graph] A verifier-accepted statement with proof; the DAG of
  facts and their dependencies (Danus, Pharos).
  \item[FNR, TNR] False-negative rate (valid solutions rejected) and
  true-negative rate (invalid solutions rejected), with valid solutions as the
  positive class.
  \item[Harness] Everything in an agent other than the model weights.
  \item[MCP] Model Context Protocol, a standard way of exposing tools to an
  agent.
  \item[Reasoning effort] A setting (\texttt{high}, \texttt{xhigh},
  \texttt{max}, \ldots) controlling how long the model thinks per step.
  \item[Route board] Pharos's record of approaches, their status, evidence,
  and obstructions.
  \item[Skill] A Markdown instruction file for a recurring task, loaded by the
  agent when relevant.
  \item[Subagent] An agent launched by another agent, with its own context.
  \item[Test-time compute] Computation spent when solving a task, as opposed
  to training.
\end{description}

\clearpage
\printbibliography

\end{document}